\input amstex
\documentstyle{amsppt}
\document

\magnification 1100

\def\?{{\bold \,(?????)}}
\def\1b{{\bold 1}}
\def\ab{{\bold a}}
\def\bb{{\bold b}}
\def\cb{{\bold c}}
\def\db{{\bold d}}
\def\eb{{\bold e}}
\def\fb{{\bold f}}
\def\hb{{\bold h}}
\def\ib{{\bold i}}
\def\jb{{\bold j}}
\def\kb{{\bold k}}
\def\mb{{\bold m}}
\def\nb{{\bold n}}
\def\lb{{\bold l}}
\def\ob{{\bold o}}

\def\Ab{{\bold A}}
\def\Bb{{\bold B}}

\def\Ib{{\bold I}}

\def\Mb{{\bold M}}
\def\Kb{{\bold K}}
\def\Hb{{\bold H}}

\def\Rb{{\bold R}}
\def\Sb{{\bold S}}
\def\Tb{{\bold T}}

\def\Ub{{\bold U}}

\def\x{{\roman x}}

\def\Hom{\text{Hom}\,}

\def\fock{\bigwedge^{\infty}}

\def\Ker{\text{Ker}\,}
\def\Im{\text{Im}\,}

\def\GL{{\text{GL}}}

\def\mod{\text{mod}\,}

\def\AA{{\Bbb A}}

\def\CC{{\Bbb C}}
\def\FF{{\Bbb F}}

\def\LL{{\Bbb L}}

\def\NN{{\Bbb N}}

\def\QQ{{\Bbb Q}}

\def\SS{{\Bbb S}}

\def\ZZ{{\Bbb Z}}

\def\Dc{{\Cal D}}

\def\Fc{{\Cal F}}
\def\Gc{{\Cal G}}
\def\Hc{{\Cal H}}

\def\Sc{{\Cal S}}

\def\Vo{{\bar V}}

\def\abo{{\bar\ab}}
\def\bbo{{\bar\bb}}
\def\dbo{{\bar\db}}
\def\ao{{\bar a}}
\def\bo{{\bar b}}
\def\do{{\bar d}}
\def\io{{\bar{\imath}}}

\def\jo{{\bar{\jmath}}}

\def\ko{{\bar k}}

\def\ts{\textstyle}                
                
\def\qed{\hfill $\sqcap \hskip-6.5pt \sqcup$}        % White box                
\overfullrule=0pt                                    % No black boxes

\def\slinf{{{\frak{sl}}_{\infty}}}
\def\slh{{\widehat{\frak{sl}}_n}}

\def\lra{{\longrightarrow}}

\def\var{\varepsilon}

\def\lra{{\longrightarrow}}
\def\lla{{\longleftarrow}}

\newdimen\Squaresize\Squaresize=14pt
\newdimen\Thickness\Thickness=0.5pt
\def\Square#1{\hbox{\vrule width\Thickness
	      \vbox to \Squaresize{\hrule height \Thickness\vss
	      \hbox to \Squaresize{\hss#1\hss}
	      \vss\hrule height\Thickness}
	      \unskip\vrule width \Thickness}
	      \kern-\Thickness}
\def\Vsquare#1{\vbox{\Square{$#1$}}\kern-\Thickness}
\def\young#1{
\vbox{\smallskip\offinterlineskip
\halign{&\Vsquare{##}\cr #1}}}

%%%%%%%%%%%%%%%%%%%%%%%%%%%%%%%%%%%%%%%%%%%%%%%%%%%%%%%%%%%%%%%%%%%%%%%%%%%%%%%%
\centerline{\bf ON THE DECOMPOSITION MATRICES}
\centerline{\bf OF THE QUANTIZED SCHUR ALGEBRA}

\vskip15mm

\centerline{\bf Michela VARAGNOLO and Eric VASSEROT
\footnote{Both authors are partially supported by the EEC grant
no. ERB FMRX-CT97-0100.}}

\vskip15mm

\noindent{\bf Abstract.} We prove the decomposition conjecture for the Schur 
algebra stated in [LT]. We also give a new approach to the Lusztig 
conjecture via canonical bases of the Hall algebra.

\vskip15mm

\centerline{\bf 0. Introduction and general notations.}

\vskip3mm

\noindent {\bf 0.1.}
The aim of this paper is to give a proof of the decomposition
conjecture for the quantized Schur algebra [LT, Conjecture 5.2]
which generalizes the theorem of Ariki (see [A])
on the decomposition numbers of the Hecke algebra of type $A$.
More precisely, let $\bigwedge^\infty$ be the level 1 Fock space of 
type $A$ and let $\Bb^\pm$ be the bases of $\bigwedge^\infty$ introduced 
in [LT]. The decomposition
conjecture links the decomposition matrices of the quantized 
Schur algebra and the basis $\Bb^+$. Our proof consists in two steps :
first we express $\Bb^\pm$ in terms of some Kazhdan-Lusztig polynomials.
Then we note that a simple module of the quantized Schur algebra 
can be pulled-back to a simple module of the Lusztig integral
form of the quantized enveloping
algebra of ${\frak{sl}}_k$ (denoted by $\Ub({\frak{sl}}_k)$).
Thus, the Lusztig conjecture for the dimension of the simple 
$\Ub({\frak{sl}}_k)$-modules at roots of unity 
identifies the entries of the decomposition matrices with some
Kazhdan-Lusztig polynomials. 
It suffices to observe that these polynomials are precisely 
the ones which appear in $\Bb^+$. 

Let $\Ub_n^-$ be the Hall algebra of nilpotent representations
of the cyclic quiver. Set $\varepsilon=\exp(2i\pi/n')$.
Put $n=n'$ if $n$ odd, $n=n'/2$ else, i.e.
$n$ is the order of $\varepsilon^2$.
Let $\Ub_\varepsilon({\frak{sl}}_k)$ be the specialization at $v=\var$ of
$\Ub(\frak{sl}_k)$.
We give a new approach to the proof of 
the Lusztig conjecture on the character of the simple modules
of $\Ub_\varepsilon({\frak{sl}}_k)$
in terms of the canonical basis of $\Ub_n^-$.
Recall that this conjecture 
(proved by Kashiwara-Tanisaki and Kazhdan-Lusztig)
gives the multiplicity
of the Weyl module of $\Ub_\varepsilon({\frak{sl}}_k)$
with highest weight $\mu$, say $W_\mu$,
in the simple $\Ub_\varepsilon({\frak{sl}}_k)$-module
with highest weight $\lambda$, say $V_\lambda$, i.e.
$$[V_\lambda:W_\mu]=\sum_y(-1)^{l(yx)}P_{yx}(1),\leqno(a)$$
where $x\in\widehat{\frak S}_k$ is minimal such that 
$\nu=\lambda\cdot x^{-1}$ satisfies
$$\nu_i<\nu_{i+1}\quad\forall i=1,2,...,k-1,\qquad\nu_1-\nu_k\geq 1-k-n,$$
and $\mu=\lambda\cdot x^{-1}y.$  
We proceed as follows.
First we prove that $\bigwedge^\infty$ is a cyclic $\Ub^-_n$-module
generated by the vacuum vector $|\emptyset\rangle.$ 
Then we define a basis $\Bb'$ of $\Ub_n^-$ using intersection cohomology.
We construct a basis $\Bb$ of $\bigwedge^\infty$ via the action of 
$\Bb'$ on the vacuum vector. 
We prove that $\Bb$ and $\Bb^+$ are fixed by the same semi-linear 
involution (see Theorem 6.3). At last, we prove that the equality $\Bb=\Bb^+$ 
is a $q$-analogue of the Lusztig conjecture (see Subsection 11.4).
The reader should be warned that we endow the Hall algebra with the 
product opposit to the usual one (used in [G1] or [L1-4]).

The plan of the paper is the following. In Sections 1-4 we recall the 
definitions and the main properties of the basic objects. In Sections 5-6
we construct an action of $\Ub^-_n$
on the Fock space $\bigwedge^\infty$. Proposition 6.1 is new. 
In Section 7 we introduce the 
convolution algebra on pairs of affine flags. This algebra is a geometric
analogue of the affine Schur algebra (Proposition 7.4) and is related to
$\Ub^-_n$ in Proposition 7.6. In Sections 8-9 we give a representation
of $\Ub_n^-$ on the finite wedges space, $\bigwedge^l$, via the coproduct
of $\Ub_n^-$. This action is related to the convolution algebra
on affine flags by Lemma 8.3. 
In Section 10 we interpret the action of
$\Ub_n^-$ on $\bigwedge^\infty$ as a ``limit" of $\bigwedge^l$ when $l$
goes to infinity. Using the results of Sections 7-9 we prove that
the elements of $\Bb$ are fixed by the Leclerc-Thibon involution
(Theorem 6.3).
In Section 11 we prove the Decomposition Conjecture.
Let us observe that the proof only uses the results of Sections 8 and 9.
In Section 12 we reinterpret the Lusztig conjecture. 
We use in an essentiel way the construction of the representation of 
$\Ub^-_n$ on $\bigwedge^\infty$ given in Section 6. 

\vskip3mm

\centerline{\bf Contents.}
\itemitem
\hfill\break
{\bf 1. The Hecke algebra. }  
\hfill\break
{\bf 2. The quantum group. }
\hfill\break
{\bf 3. The Hall algebra.}
\hfill\break
{\bf 4. The Fock space.} 
\hfill\break
{\bf 5. The representation of $\Ub^-_\infty$ on $\fock$.}
\hfill\break
{\bf 6. The representation of $\Ub^-_n$ on $\fock$.}
\hfill\break
{\bf 7. Flag varieties.}
\hfill\break
{\bf 8. The tensor representation of $\tilde\Ub^-_n$.}
\hfill\break
{\bf 9. The action of $\Ub^-_n$ on wedges.}
\hfill\break
{\bf 10. Proof of Theorem 6.3.}
\hfill\break
{\bf 11. Proof of the Decomposition Conjecture.}
\hfill\break
{\bf 12. The Lusztig conjecture.}
\hfill\break
{\bf 13. Proof of Proposition 6.1.}

\vskip3mm

\noindent{\bf 0.2.}
We now fix a few general notations. Set $\SS=\CC[v]$, $\AA=\CC[v,v^{-1}].$ 
Let $\FF$ be a field with $q^2$ elements and let $\bar\FF$ be the
algebraic closure of $\FF.$
Fix a set $I$. For any $i\in I$ and $r\in\NN^\times$,
let $\bar\FF^r[i]$ be the $I$-graded $\bar\FF$-vector space with a single
$r$-dimensional component, in degree $i$.
Let $\epsilon_i\in\NN^{(I)}$ be the dimension of $\bar\FF[i]$.
For any $d\in\NN^{(I)}$ set $|d|=\sum_{i\in I}d_i$.
If $i\in\ZZ$ let $\io$ be the class of $i$ in $\ZZ/n\ZZ$. 
Given a positive integer $l$ let $\Pi(l)$ be the set of all the
partitions of $l$ and let $\Pi_l$ be the set of partitions with
at most $l$ parts. Put $\Pi=\cup_l\Pi(l)$.
The set $\Pi$ is endowed with the usual order.
If $\lambda\in\Pi$ let $\lambda'$ be the dual partition.
For an irreducible algebraic variety $X$ we denote by $\Hc^i(IC_X)$ the
$i$-th cohomology sheaf of the intersection complex of $X$.
Then, for any stratum $Y\subset X$, let $\dim\Hc^i_Y(IC_X)$ be the dimension
of the stalk of $\Hc^i(IC_X)$ at a point of $Y$.
For any set $X$ with the action of a group $G$
let $\CC_G(X)$ be the set of $G$-invariant functions $X\to\CC$
supported on a finite number of orbits.
For any subset $X$ of an algebraic variety let $\bar X$ denote its
Zariski closure. 

\vskip5mm

\centerline{\bf 1. The Hecke algebra. }  

\vskip3mm

\noindent {\bf 1.1.}
Fix $n\in\NN^\times$ and set
$$A_l^n=\{\ib\in\ZZ^l\,|\,1-n\leq i_1\leq i_2\leq\cdots\leq i_l\leq 0\}.$$
Let ${\frak S}_l$ be the symmetric group and let 
$\widehat{\frak S}_l={\frak S}_l\ltimes\ZZ^l$ be the extended
affine Weyl group. Let $\widehat S_l\subset\widehat{\frak S}_l$ 
be the set of simple affine reflexions and put
$S_l=\widehat S_l\cap{\frak S}_l$. 
As usual, the simple affine reflexions are denoted by $s_0,s_1,...,s_{l-1}$ 
in such a way that $S_l=\{s_1,s_2,...,s_{l-1}\}$.
Let $\pi\in\widehat{\frak S}_l$ be the zero length element such that
$s_{i-1}=\pi^{-1}s_i\pi.$
The group $\widehat{\frak S}_l$ acts on $\ZZ^l$ on the right in such a way
that  
$$\matrix
(\ib)\lambda=\ib+n\lambda\qquad\hfill&\text{if}\quad\lambda\in\ZZ^l\hfill\cr\cr
(\ib)s_j=(i_1,i_2,...,i_{j+1},i_j,...,i_l)\qquad\hfill&\text{if}\quad
j\neq 0\hfill\cr\cr
(\ib)s_0=(i_l-n,i_2,...,i_{l-1},i_1+n).\hfill&
\endmatrix $$
The alcove $A_l^n$ is a fundamental domain for this action.
If $\ib\in A_l^n$ let ${\frak S}_\ib\subset\widehat{\frak S}_l$
be its isotropy group, $S_\ib=\widehat S_l\cap{\frak S}_\ib$,
and let ${\frak S}^\ib$ 
be the set of minimal length representatives of the cosets in 
${\frak S}_\ib\setminus\widehat{\frak S}_l.$ For any $x\in\widehat{\frak S}_l,$
let $x_\ib\in{\frak S}_\ib$ and $x^\ib\in{\frak S}^\ib$ be such that 
$x=x_\ib x^\ib.$ Let $\omega\in{\frak S}_l$ be the longest element. 
Set $\rho=(0,-1,-2,...,1-l)\in\ZZ^l$ and put
$$\lambda\cdot x=(\lambda+\rho)x-\rho,\qquad x\in\widehat{\frak S}_l,\quad
\forall\lambda\in\ZZ^l.$$

\vskip1mm

\noindent{\bf 1.2.}
The Hecke algebra of type $GL_l$, say $\Hb_l$, is the 
unital associative $\AA$-algebra generated by $T_i^{\pm 1}$, $i=1,2,...l-1$
modulo the following relations
$$\matrix
T_i\,T_i^{-1}=1=T_i^{-1}\,T_i,\qquad (T_i+1)(T_i-v^{-2})=0,\cr\cr
T_i\,T_{i+1}\,T_i=T_{i+1}\,T_i\,T_{i+1},\qquad
|i-j|>1\Rightarrow T_i\,T_j=T_j\,T_i.\endmatrix\leqno(a)$$

\noindent The affine Hecke algebra of type $GL_l$, say $\widehat\Hb_l$,
is the unital associative $\AA$-algebra generated by $T_i^{\pm 1},X_j^{\pm 1},$
$i=1,2,...,l-1,j=1,2,...,l$ modulo the relations $(a)$ and
$$X_i,X_i^{-1}=1=X_i^{-1}\,X_i,\qquad X_i\,X_j=X_j\,X_i,$$
$$T_i\,X_i\,T_i=v^{-2}X_{i+1},\qquad 
j\not= i,i+1\Rightarrow X_j\,T_i=T_i\,X_j.$$
For all $x
\in{\frak S}_l\ltimes\ZZ^l$ let $l(x)$ be the length of $x$ and 
let $\tilde T_x$ be the normalized element $\tilde T_x=v^{l(x)}T_x$.
The algebra $\widehat\Hb_l$ is isomorphic to the Hecke algebra of the 
extended affine Weyl group ${\frak S}_l\ltimes\ZZ^l$ via the
Bernstein isomorphism which maps $\tilde T_\lambda^{-1}$ to
$X^\lambda=X_1^{\lambda_1}X_2^{\lambda_2}\cdots X_l^{\lambda_l}$ 
if $\lambda\in\ZZ^l$ is dominant, i.e. if
$\lambda_1\geq \lambda_2\geq \cdots\geq \lambda_l$.
The semilinear involution 
$\ \bar{}\ :\widehat\Hb_l\to\widehat\Hb_l$
is such that $\bar T_x=T^{-1}_{x^{-1}}$ for all $x$.
For all $x$ put $\tilde T_x=v^{l(x)}T_x$.
If $t\in\CC^\times$ let $\widehat\Hb_{l|t}$
be the specialization of $\widehat\Hb_l$ at $v=t$.

\vskip5mm

\centerline{\bf 2. The quantum group. }

\vskip3mm

\noindent
Put $I=\{1,2,...,n-1\}$ (resp. $I=\{0,1,...,n-1\}$) and let
$a_{ij}$ be the entries of the Cartan matrix of type $A_{n-1}$ 
(resp. $A_{n-1}^{(1)}$).
The quantized enveloping enveloping algebra of $\frak{sl}_n$ (resp.
$\widehat{\frak{sl}}_n$) is the unital associative $\CC(v)$-algebra generated
by $\eb_i,\fb_i,\kb_i^{\pm 1}$, $i\in I$, modulo the Kac-Moody type
relations
$$\kb_i\,\kb_i^{-1}=1=\kb_i^{-1}\,\kb_i,\qquad \kb_i\,\kb_j=\kb_j\,\kb_i,$$
$$\kb_i\,\eb_j=v^{a_{ij}}\eb_j\,\kb_i,\quad
\kb_i\,\fb_j=v^{-a_{ij}}\fb_j\,\kb_i,\quad
[\eb_i,\fb_j]=\delta_{ij}{\kb_i-\kb_i^{-1}\over v-v^{-1}},$$
$$\sum_{k=0}^{1-a_{ij}}(-1)^{k}\eb_i^{(k)}\,\eb_j\,\eb_i^{(1-a_{ij}-k)}=
\sum_{k=0}^{1-a_{ij}}(-1)^{k}\fb_i^{(k)}\,\fb_j\,\fb_i^{(1-a_{ij}-k)}=0
\quad\text{if}\quad i\neq j,$$
where 
$$[k]={v^k-v^{-k}\over v-v^{-1}},\qquad
[k]!=[k]\,[k-1]\,\cdots [1],\qquad
\eb_i^{(k)}={\eb_i^k\over [k]!},\qquad
\fb_i^{(k)}={\fb_i^k\over [k]!}.$$
We denote by $\Ub({\frak{sl}}_n)$ (resp. $\Ub(\widehat{\frak{sl}}_n)$) the 
Lusztig integral form, i.e. the $\AA$-subalgebra
generated by the divided powers $\eb_i^{(k)}$, $\fb_i^{(k)}$,
and by $\kb_i^{\pm 1}$. 
If $n=\infty$ the algebra $\Ub({\frak{sl}}_\infty)$ is well defined.
The algebras above are Hopf algebras. The coproduct is
$$\Delta\eb_i=\eb_i\otimes\kb_i+1\otimes\eb_i,\qquad
\Delta\fb_i=\fb_i\otimes 1+\kb_i^{-1}\otimes \fb_i,\qquad
\Delta\kb_i=\kb_i\otimes\kb_i.$$
Let $\Ub^-(\widehat{\frak{sl}}_n)\subset\Ub(\widehat{\frak{sl}}_n)$
and $\Ub^-({\frak{sl}}_\infty)\subset\Ub({\frak{sl}}_\infty)$ be the
subalgebras generated by the elements $\fb_i^{(k)}$.

\vskip5mm

\centerline{\bf 3. The Hall algebra.}

\vskip3mm

\noindent In this section we recall some of the results of [L1-4] and [G1].

\vskip3mm

\noindent{\bf 3.1.}
Fix a finite field $\FF$ with $q^2$ elements as in the introduction.
Let $\Gamma=(I,J)$ be an oriented graph : 
$I$ is the set of vertices and $J$ is the set of arrows. 
Given an arrow $j\in J$ let $j_1$ and $j_2$ be respectively
the input vertex and the output vertex. Fix $d\in\NN^{(I)}$ and let $V$
be an $I$-graded $\FF$-vector space of dimension $d$.
Let $E_V\subseteq\bigoplus_{j\in J}\Hom(V_{j_1}, V_{j_2})$
be the subset of nilpotent representations of $\Gamma$ on $V$.
In this paper we will suppose that $\Gamma$ is one of the 
following two graphs :

\vskip2mm

\noindent$(a)\quad\Gamma=\Gamma_n$ is the cyclic quiver of type $A_n^{(1)}$,
i.e. $I=\ZZ/n\ZZ$ and $J=\{\io\to\io+1\,|\,\io\in\ZZ/n\ZZ\}$,

\vskip2mm

\noindent$(b)\quad\Gamma=\Gamma_\infty$ is the infinite quiver of type 
$A_\infty$, i.e. $I=\ZZ$ and $J=\{i\to i+1\,|\,i\in\ZZ\}$.

\vskip2mm

\noindent{\bf 3.2.}
Set $\Ab_d=\CC_{G_V}(E_V)$ where $G_V=\prod_{i\in I} GL(V_i)$.
Given $a,b\in\NN^{(I)}$ such that $d=a+b$, fix $I$-graded
$\FF$-vector spaces $U,W$ of dimensions $a,b$. 
Let consider the diagram 
$$E_{U}\times E_{W}{\buildrel p_1\over\lla}
E{\buildrel p_2\over\lra}F{\buildrel p_3\over\lra}E_V,$$
where

\vskip2mm

\noindent{$(c)$} $E$ is the set of triples $(x,\phi,\psi)$
such that $x\in E_V$,
$$0\to U{\buildrel\phi\over\lra}V{\buildrel\psi\over\lra}W\to 0$$
is an exact sequence of $I$-graded vector spaces and
$\phi(U)$ is stable by $x$, 

\vskip2mm

\noindent{$(d)$} $F$ is the set of pairs $(x,U')$ where $x\in E_V$ and 
$U'\subset V$ is a $x$-stable $I$-graded subspace of dimension $a$.

\vskip2mm

\noindent Given $f\in\CC_{G_{U}}(E_{U})$ and 
$g\in\CC_{G_{W}}(E_{W})$ set 
$$f\circ g=q^{-m(b,a)}(p_3)_!h\in\CC_{G_V}(E_V),$$ 
where $h\in\CC(F)$ is the function such that $p_2^*h=p_1^*(fg)$ and
$m(b,a)=\sum_{j\in J}b_{j_1}a_{j_2}+\sum_{i\in I}b_ia_i.$
Then $(\Ab,\circ)$, where $\Ab=\bigoplus_d\Ab_d$, is an associative algebra. 

\vskip2mm

\noindent{\bf 3.3.}
Given $a,b\in\NN^{(J)}$ such that $d=a+b$, fix a $I$-graded
$\FF$-vector space $U\subset V$ of dimension $a$. 
Let consider the diagram
$$E_{U}\times E_{V/U}{\buildrel p\over\lla}
E{\buildrel i\over\lra}E_V.$$
Here $E\subset E_V$ is the subset of representations
preserving $U$, the map $i$ is the inclusion and
$p$ is the obvious projection. Set
$$\Delta_{a,b}\,:\,\Ab_d\to\Ab_a\otimes\Ab_b,\quad
f\mapsto q^{-n(b,a)}p_!i^*f,$$
where $n(b,a)=\sum_{j\in J}b_{j_1}a_{j_2}-\sum_{i\in I}b_ia_i$.

\vskip3mm

\noindent{\bf 3.4.}
Recall that $\Gamma=\Gamma_n$ or $\Gamma_\infty$.
The classification of the isomorphism classes of nilpotent 
representations of $\Gamma$ does not depend on  the ground field $\FF.$
It is proved in [R] that the structural constants of $\Ab$ in the basis
formed by the characteristic functions of the $G_V$-orbits in $E_V$ 
are the value at $v=q$ of universal polynomials in $\AA.$
Thus $\Ab$ can be viewed as the specialization at $v=q$ 
of a $\AA$-algebra, called the generic Hall algebra. 
Let $\Ub^-_n$ (resp. $\Ub^-_\infty$) be the generic
Hall algebra if  $\Gamma=\Gamma_n$ (resp. $\Gamma=\Gamma_\infty$).
It is known that $\Ub^-_\infty$ is isomorphic to $\Ub^-(\slinf)$
and that $\Ub^-(\slh)$ embeds in $\Ub^-_n$ (see [G1]).
Let $\Ab^0$ be the $\AA$-linear span of elements
$\kb_d$ with $d\in\ZZ^{(I)}$ such that
$$\kb_0=1\quad\text{and}\quad
\kb_{a}\kb_{b}=\kb_{a+b},\quad\forall a,b.$$
For simplicity we will write $\kb_i=\kb_{\epsilon_i}$ for all $i\in I$.
Set $\tilde\Ab=\Ab\otimes_\AA\Ab^0$ and put
$$(f\otimes\kb_a)\circ(g\otimes\kb_b)=v^{-a\cdot d}
(f\circ g)\otimes\kb_{a+b},\qquad\forall g\in\Ab_{d}\quad\forall f\in\Ab,$$ 
where $a\cdot d=-n(a,d)-n(d,a).$ Consider the map
$\Delta\,:\,\tilde\Ab\to\tilde\Ab\otimes_\AA\tilde\Ab$
such that
$$\Delta(f\otimes\kb_c)=\sum_{d=a+b}
\Delta_{a,b}(f)(\kb_{b+c}\otimes\kb_c),\quad\forall f\in\Ab_d.$$
Then $(\tilde\Ab,\circ,\Delta)$ is a $\AA$-bialgebra 
(it is due to Lusztig for the composition algebra, 
the general case is due to Green).
Put $\tilde\Ub^-_n=\tilde\Ab$ if $\Gamma=\Gamma_n$ 
and $\tilde\Ub^-_\infty=\tilde\Ab$ if $\Gamma=\Gamma_\infty$.

\vskip3mm

\noindent{\bf 3.5.}
Given a $G_V$-orbit $O\subset E_V$ let $\fb_O\in\Ab$ 
be the $v^{\dim\,O}$ times the characteristic function of $O$. 
For any $G_V$-orbit $O\subset E_V$ set 
$$\bb_O=\sum_{i,O'}v^{-i+\dim O-\dim O'}\dim\Hc^i_{O'}(IC_O)\,\fb_{O'}.$$
The elements $\bb_O$ form a basis of $\Ab$.
If $d\in\NN^{(I)}$ let $\fb_d\in\Ab$ be the characteristic function 
of the zero representation of $\Gamma$ in a $d$-dimensional space. 
The following result is proved in Section 13.

\vskip3mm

\noindent{\bf Proposition.} {\it The algebra $\Ab$ is generated by the $\fb_d$,
$d\in\NN^{(I)}$.}
\qed

\vskip3mm

\noindent{\bf 3.6.}
Given two integers $i\leq j$, let $\bar\FF[i,j]$ be the unique indecomposable
representation of $\Gamma_\infty$ (resp. $\Gamma_n$)
with dimension $\sum_{k=i}^j\epsilon_k$ (resp. $\sum_{k=i}^j\epsilon_\ko$).
For any partition $\lambda=(\lambda_1\geq\lambda_2\geq\cdots)$ 
let $\bar\FF[\lambda]$ be the representation of $\Gamma$ such that
$$\bar\FF[\lambda]=\bigoplus_{k\geq 1}\bar\FF[1-k,\lambda_k-k].$$
Let $O_\lambda$ be the orbit of $\bar\FF[\lambda]$ and put 
$d_\lambda=\dim O_\lambda$.

\vskip5mm

\centerline{\bf 4. The Fock space.} 

\vskip3mm

\noindent In this section we recall the construction of the 
quantized Fock space, due to [H], as it is re-interpreted in [MM].

\vskip3mm

\noindent{\bf 4.1.}
Let $T(\lambda)$ be the tableau of shape $\lambda$ whose box with 
coordinates $(x,y)$ is filled with $y-x$. For instance if 
$\lambda=(432)$ we get

$$\young{-2&-1\cr -1&0&1\cr 0&1&2&3\cr}.$$

\noindent Let $\fock$ be a $\AA$-module with basis
$\{|\lambda\rangle\,|\,\lambda\in\Pi\}$.
If $i\in\ZZ$, a removable $i$-box of $T(\lambda)$ is a box with the color $i$
which can be removed in such a way that the
new tableau still comes from a partition.
Similarly, an indent $i$-box corresponds to a box with the color $i$
which can be added to $T(\lambda)$.
Given $\io\in\ZZ/n\ZZ$, $i\in\io$, and a partition $\lambda$ put

$$n_i(\lambda)=\sharp\{\text{indent\ $i$-box\  of\ }T(\lambda)
\}-\sharp\{\text{removable\ $i$-box\ of\ }T(\lambda)\},$$

\noindent and
$n_\io(\lambda)=\sum_{i\in\io}n_i(\lambda),$
$n^-_i(\lambda)=\sum_{j<i\&j\in\io}n_j(\lambda),$
$n^+_i(\lambda)=\sum_{j>i\&j\in\io}n_j(\lambda).$

\vskip3mm

\noindent{\bf 4.2.}
The algebra $\Ub(\slinf)$ acts on $\fock$ by
$$\kb_i(|\lambda\rangle)=v^{n_i(\lambda)}\,|\lambda\rangle,
\quad \eb_i(|\lambda\rangle)=|\nu\rangle,
\quad \fb_i(|\lambda\rangle)=|\mu\rangle,$$
where the partitions $\mu,\nu$ are 
such that $T(\mu)-T(\lambda)$ and $T(\lambda)-T(\nu)$
are a box with color $i$. It is known that
$\fock$ is the simple module with highest weight $\Lambda_0$ 
(the fundamental weight) and that the canonical basis of $\fock$
is $\{|\lambda\rangle\,|\,\lambda\in\Pi\}.$
The weight multiplicities in $\fock$ are $0$ or $1$,
i.e. $\Lambda_0$ is a minuscule weight. 

\vskip3mm

\noindent{\bf 4.3.}
The algebra $\Ub(\slh)$ acts on $\fock$ by 
$$\kb_\io(|\lambda\rangle)=v^{n_\io(\lambda)}\,|\lambda\rangle,
\quad
\eb_\io(|\lambda\rangle)=
\sum_{i\in\io}v^{-n^-_i(\lambda)}\eb_i(|\lambda\rangle),
\quad
\fb_\io(|\lambda\rangle)=
\sum_{i\in\io}v^{n^+_i(\lambda)}\fb_i(|\lambda\rangle).$$

\vskip5mm

\centerline{\bf 5. The representation of $\Ub^-_\infty$ on $\fock$.}

\vskip3mm

\noindent
The algebras $\Ub^-_\infty$ and $\Ub^-(\slinf)$ are isomorphic.
Thus $\fock$ may be viewed as the quotient of 
$\Ub^-_\infty$ by a left ideal $\Ib$. Let us describe $\Ib$.
Let $\bar\Gamma_\infty$ be the quiver $\Gamma_\infty$
with the opposit orientation. For any $\ZZ$-graded 
$\bar\FF$-vector space $V$ let $\Lambda_V$ be the variety of pairs $(x,\bar x)$ 
of commuting representations respectively of $\Gamma_\infty$ and 
$\bar\Gamma_\infty$ on $V$. The variety $\Lambda_V$ is reducible.
For any $G_V$-orbit $O\subset E_V$ set
$$\Lambda_O=\{(x,\bar x)\in\Lambda_V\,|\,x\in O\}.$$
According to [N] the orbit $O$ is stable if there exists a triple 
$$(x,\bar x,i)\in\overline{\Lambda}_O\times\text{Hom}(\bar\FF[0],V)$$
such that $i$ is homogeneous of degree 0 and that
for any graded subspace $W\subseteq V$,
$$(x(W),\bar x(W)\subseteq W \quad\text{and}\quad \Im i\subseteq W)
\quad\Rightarrow\quad W=V\leqno(a)$$
(since the Hall algebra is endowed with the product opposit to the usual one,
we use the stability condition opposit to the one in [N]).

\vskip1mm

\noindent{\bf Proposition.} {\it
The ideal $\Ib$ is linearly spanned by the elements $\bb_O$ such that
$O\neq O_\lambda$ for all $\lambda$. Moreover
the map $\Ub^-_\infty/\Ib\to\fock,\,
\bb_{O_\lambda}+\Ib\mapsto|\lambda\rangle$, is an isomorphism
of $\Ub^-_\infty$-modules.}

\vskip3mm

\noindent{\it Proof.} 
From [N, Theorem 11.7 and Proposition 3.5], $\Ib$ 
is linearly generated by the elements $\bb_O$ such that $O$ is unstable.
Let us show that for any $\lambda\in\Pi$ the orbit $O_\lambda$ is
stable. A dimension counting then shows that the
orbits $O_\lambda$ are precisely all the stable orbits.
Recall that $\bar\FF[i,j]$ is the representation 
$x$ of $\Gamma_\infty$ on the graded space $\bigoplus_{k=i}^j\bar\FF\, v_k$, 
where $v_k$ is a non-zero vector of degree $k$, such that
$x(v_k)=v_{k+1}$ if $k<j$ and $x(v_j)=0.$
Fix $\lambda=(\lambda_1,\lambda_2,...,\lambda_r)\in\Pi$.
Fix non zero vectors $v_{k,s}\in\bar\FF[1-k,\lambda_k-k]$ with degree $s$.
The representation $\bar\FF[\lambda]$ is given 
by the endomorphism $x$ such that for all $k$,
$$x(v_{k,s})=v_{k,s+1}\quad\text{if}\quad s\in [1-k,\lambda_k-k), 
\quad\text{and}\quad x(v_{k,\lambda_k-k})=0.$$
Let us exhibit a pair $(i,\bar x)$ satisfying $(a)$.
Fix a graded homomorphism $i\in\text{Hom}(\bar\FF[0],\bar\FF[\lambda])$ 
such that $v_{1,0}\in\Im i$.
Consider the degree $-1$ linear operator $\bar x$ on $\bar\FF[\lambda]$ such that
$$\bar x(v_{k,s})=v_{k+1,s-1}\ \text{if}\ k\not=r\ \text{and} \ 
s\leq\lambda_{k+1}-k,\qquad \bar x(v_{k,s})=0\ \text{else.}$$
The operators $x$ and $\bar x$ commute since
$$\matrix
x\bar x(v_{r,s})=0=\bar x x(v_{r,s})\qquad\hfill&\forall s,
\hfill\cr\cr
x\bar x (v_{k,s})=0=\bar x x(v_{k,s})\qquad\hfill&
\forall s\geq \lambda_{k+1}-k,
\hfill\cr\cr
x\bar x(v_{k,s})=v_{k+1,s}=\bar x x(v_{k,s})\qquad\hfill&
\forall s<\lambda_{k+1}-k,\ \forall k\not=r.\hfill
\endmatrix$$
Now if $W\subseteq V$ is such that $x(W)\subseteq W$ and $\Im i\subseteq W$
then $\bar\FF[0,\lambda_1-1]\subseteq W$ : namely $v_{1,0}\in W$ and thus 
$v_{1,s}=x^s(v_{1,0})\in W$ for all $s$.
By definition of $\bar x$ we have for all $t<r$
$$\bar x^t(\bar\FF[0,\lambda_1-1])=\bar\FF[-t,\lambda_{1+t}-1].$$
The dimension $d_\lambda$ of $\bar\FF[\lambda]$ is
such that $d_{\lambda,i}$ is the multiplicity of the color $i$ 
in the tableau $T(\lambda)$ (see Section 4). Thus
the linear isomorphism $\bb_{O_\lambda}\mapsto |\lambda\rangle$
preserves the weights. Moreover it preserves the canonical base up to a 
permutation of its elements. Since $\Lambda_0$ is minuscule there is at most 
one vector of a given weight in the canonical basis. Hence 
the canonical bases are fully identified.
\qed

\vskip5mm

\centerline{\bf 6. The representation of $\Ub^-_n$ on $\fock$.}

\vskip3mm

\noindent{\bf 6.1.}
Fix $d\in\NN^{(\ZZ)}$ and let
$V$ be a $\ZZ$-graded $\bar\FF$-vector space of dimension $d$. Let 
$\do\in\NN^{\ZZ/n\ZZ}$ be the multi-index such that
$\do_\io=\sum_{j\in\io}d_j$ for all $\io\in\ZZ/n\ZZ,$
and let $\Vo$ be the $\do$-dimensional $\ZZ/n\ZZ$-graded vector 
space such that $\Vo_\io=\bigoplus_{j\in\io}V_j.$
The vector space $\Vo$ is filtered by the subspaces 
$$\Vo_{\geq i}=\bigoplus_{j\geq i}V_j,\qquad\forall i\in\ZZ.$$
The associated graded is naturally identified with the $\ZZ$-graded space $V$.
Set
$$E_{\Vo,V}=\{x\in E_\Vo\,|\,x(\Vo_{\geq i})\subseteq\Vo_{\geq i+1},
\quad\forall i\}.$$
The map $p\,:\,E_{\Vo,V}\to E_V$ associate to a representation of
$\Gamma_n$ in $\Vo$ the corresponding graded representation
of $\Gamma_\infty$ in $V$. Let 
$j\,:\,E_{\Vo,V}\hookrightarrow E_{\Vo}$ be the closed embedding. 
Let consider the map $\gamma_d\,:\,\Ub^-_{n,\do}\to\Ub^-_{\infty,d}$ 
such that 
$$\gamma_{d|v=q^{-1}}\,:\,\CC_{G_\Vo}(E_\Vo)\to\CC_{G_V}(E_V),
\quad f\mapsto q^{-h(d)}p_!j^*(f),$$
where $h(d)=\sum_{i<j\&\io=\jo}d_i(d_{j+1}-d_j).$ 
Put $k(b,a)=\sum_{i>j\&\io=\jo}b_i(2a_j-a_{j-1}-a_{j+1}).$ 
The following is proved in Section 13.

\vskip3mm

\noindent{\bf Proposition.} {\it 
Fix $\alpha,\beta\in\NN^{\ZZ/n\ZZ}$ and $d\in\NN^{(\ZZ)}$ such that
$\do=\alpha+\beta$. Then,
$$\sum_{a+b=d\atop\bar a=\alpha,\bar b=\beta}
v^{-k(b,a)}\gamma_{a}(f)\circ\gamma_{b}(g)=
\gamma_{d}(f\circ g)\qquad
\forall f\in\Ub^-_{n,\alpha},\forall g\in\Ub^-_{n,\beta},$$ 
}
\qed

\vskip3mm

\noindent{\bf Remark.}
With the notations in Section 3.5 we have $\gamma_d(\fb_\do)=v^{h(d)}\fb_d$.
Observe that $\fb_d$ is the product of the divided 
powers $\fb_i^{(d_i)}$'s ordered from $i=-\infty$ to $\infty$. 

\vskip3mm

\noindent{\bf 6.2.}
For all $\lambda\in\Pi$ and all $x\in\Ub^-_n$ put 
$$x(|\lambda\rangle)=\sum_d\gamma_d(x)\kb_{d'}(|\lambda\rangle)
\quad\text{where}\quad d'=\sum_{j<i,\io=\jo}d_j\epsilon_i.\leqno(a)$$

\noindent{\bf Corollary.} {\it
Formula $(a)$ extends the Hayashi action of $\Ub^-(\slh)$ 
on $\fock$ to a representation of $\Ub^-_n$. 
}

\vskip3mm

\noindent{\it Proof.} 
The compatibility with the product in $\Ub_n^-$ follows from Proposition 6.1. 
Formula $(a)$ implies that
$$\fb_\io(|\lambda\rangle)=
\sum_{i\in\io}\sum_\mu v^{g(\epsilon_i,d_\lambda)}|\mu\rangle,$$
where $\mu$ is a partition such that $T(\mu)-T(\lambda)$ is a box with color $i$
and
$$g(\epsilon_i,d_\lambda)=-\sum_{i<j\atop\jo=\io}
(2d_{\lambda,j}-d_{\lambda,j-1}-d_{\lambda,j+1})+\alpha_i,$$
where $\alpha_i=1$ if $i<0$ and $\io=0$, and $\alpha_i=0$ else.
We have already observed that $d_{\lambda,i}$ is the multiplicity of the color 
$i$ in $T(\lambda)$. Thus,
$n_j(\lambda)=-2d_{\lambda,j}+d_{\lambda,j-1}+d_{\lambda,j+1}+\delta_{j0},$
and
$$g(\epsilon_i,d_\lambda)=\sum_{i<j\atop\jo=\io}n_j(\lambda)=n^+_i(\lambda).$$
\qed

\noindent{\bf 6.3.} 
For any $\lambda\in\Pi$ set $\bb_\lambda=\bb_{O_\lambda}|\emptyset\rangle$ 
where $O_\lambda$ is the isomorphism class of representations of $\Gamma_n$
defined in Subsection 3.6. Put $\Bb=\{\bb_\lambda\,|\,\lambda\in\Pi\}.$
Leclerc and Thibon have introduced in [LT] a semi-linear involution on 
$\fock$.

\vskip3mm

\noindent{\bf Theorem.} {\it
$\Bb$ is a basis of $\fock$ whose elements are
fixed by the Leclerc-Thibon involution.}

\vskip3mm

\noindent The theorem is proved in Subsection 10.1. 
We first introduce some more material.

\vskip3mm

\noindent{\bf 6.4.}
Let $r\,:\,E_V\to E_{\bar V}$ be such that
$${\ts r(x)_{|\Vo_\io}=\bigoplus_{i\in\io}\,x_{|V_i}\qquad\forall x\in E_V.}$$ 
Fix a $G_{\bar V}$-orbit $O\subset E_{\bar V}$
such that $O\cap E_{\bar V, V}\neq\emptyset$.
If $x\in p(O\cap E_{\bar V,V})$ then 
$$\sharp(p^{-1}(x)\cap O)\in 
\left\{\matrix
q^{2\NN}\quad\hfill&\text{if}\quad r(x)\in O\hfill\cr
(q^2-1)\NN\quad\hfill&\text{else}.\hfill
\endmatrix\right.\leqno(b)$$ 
Indeed, fix $y\in p^{-1}(x)\cap O$.
It suffices to consider the case where $y$ is indecomposable.
Then, fix a basis of homogeneous vectors $\{v_i\,|\,i\in[1,r]\}$
of $\bar V$ such that
$$y(v_k)=v_{k+1}\qquad\forall k=1,2,...,r-1.\leqno(c)$$
If $r(x)\in\bar O\setminus O$ then there exist $i,j,$ such that
$v_i\in\bar V_{\geq j}\setminus\bar V_{>j}$ and
$v_{i+1}\in\bar V_{\geq j}$.
If $t\in\FF^\times$ the representation $y_t\in E_{\bar V,V}$ 
obtained by doing $v_k\mapsto tv_k$ for all $k\leq i$ in $(c)$ 
is in $p^{-1}(x)\cap O$.
Thus $\sharp(\FF^\times)\,|\,\sharp(p^{-1}(x)\cap O).$ 
If $r(x)\in O$ then $r(x)$ and $y$ are isomorphic since
$r(x),y\in O$. Thus $x$ is indecomposable and 
it is easy to see that $p^{-1}(x)\cap O$ is a vector space. We are done. 
The identity $(b)$ implies the following lemma which is used in Section 12. 

\vskip3mm

\noindent{\bf Lemma.} {\it
For any $G_\Vo$-orbit $O\subseteq E_\Vo$ we have 
$$\gamma_d(\fb_O)=\fb_{O'}\mod\, (v-1)$$
where $O'\subseteq E_V$ is the unique $G_V$-orbit 
such that $r(O')\subseteq O$.}\qed

\vskip5mm

\centerline{\bf 7. Flag varieties.}

\vskip3mm

\noindent {\bf 7.1.}
Fix a positive integer $l$. Set $\LL=\FF((z))$ and $G=GL_l(\LL)$.
A lattice in $\LL^l$ is a free $\FF[[z]]$-submodule of rank $l$.
Let $Y$ be set of sequences of lattices $L=(L_i)_{i\in\ZZ}$ such that
$$L_i\subseteq L_{i+1}\quad\text{and}
\quad L_{i+n}=z^{-1}\,L_i.$$
The group $G$ acts on $Y$ in the obvious way.
Let $M$ be the set of all $\ZZ\times\ZZ$-matrices with non-negative
entries, say $\mb=(m_{ij})_{i,j\in\ZZ}$, such that $m_{i+n,j+n}=m_{ij}$. Set
$${\ts M^l=\{\mb\in M\,|\,\sum_{i\in\ZZ}\sum_{j=1}^nm_{ij}=l\}.}$$
The set $M^l$ parametrizes the orbits of the diagonal action of
$G$ in $Y\times Y$ : to $\mb$ corresponds the set $Y_\mb$ of
the pairs $(L',L)$ such that
$$m_{ij}=\dim_{\FF}\biggl({L_{i+1}\cap L'_{j+1}\over
(L_i\cap L'_{j+1})+(L_{i+1}\cap L'_j)}\biggr).$$
For all $L\in Y$ let $Y_{\mb,L}$ be the fiber over $L$ of the first 
projection $Y_\mb\to Y$. If $Y_{\mb,L}\neq\emptyset$ then
$Y_{\mb,L}$ is the set of $\FF$-points of an algebraic variety 
whose dimension, denoted by $y(\mb)$, is independent of $L$.
Let $\1b_\mb\in\CC_G(Y\times Y)$ be 
$q^{-y(\mb)}$ times the characteristic function of $Y_\mb$.
The convolution product, denoted $\star$, endows $\CC_{G}(Y\times Y)$
with the structure of an associative algebra. 

\vskip3mm

\noindent{\bf 7.2.}
Let $X$ be the set of sequences of lattices $L=(L_i)_{i\in\ZZ}$ such that
$$L_i\subseteq L_{i+1},\quad L_{i+l}=z^{-1}\,L_i
\quad\text{and}\quad \dim_\FF(L_{i+1}/L_i)=1.$$
The group $G$ acts on $X$ in the obvious way.
The orbits of the diagonal action of $G$ in $Y\times X$ are labelled
by functions $\ib\,:\,\ZZ\to\ZZ$ such that $\ib(k+l)=\ib(k)+n$ for all $k$ : 
let $X_\ib$ be the orbit of the pair $(L_\ib,L_\emptyset)$ such that
$$L_{\ib,i}=\prod_{\ib(j)\leq i}\FF\, e_j
\quad\text{and}\quad
L_{\emptyset,i}=\prod_{j\leq i}\FF\, e_j.$$
Here $(e_1,e_2,...,e_l)$ is a fixed $\LL$-basis of $\LL^l$ 
and $e_{i+lk}=z^{-k}\, e_i$ for all $k\in\ZZ$.
A periodic function $\ib$ as above is identified with the $l$-uple
$(\ib(1),\ib(2),...,\ib(l))\in\ZZ^l$.
If $L\in Y$ let $X_{\ib,L}$ be the fiber over $L$ of the 
projection $X_\ib\to Y$. If $X_{\ib,L}\neq\emptyset$, then $X_{\ib,L}$ 
is the set of $\FF$-points of an algebraic variety of dimension
$l(\omega_\ib)$. Let $\1b_\ib\in\CC_G(Y\times X)$ be 
$q^{-l(\omega_\ib)}$ times the characteristic function of $X_\ib$.
The space $\CC_{G}(Y\times X)$ is a left $\CC_{G}(Y\times Y)$-module
and a right $\CC_{G}(X\times X)$-module.

\vskip3mm

\noindent{\bf 7.3.}
For all $x\in\widehat{\frak S}_l$ let $X_x\subset X\times X$ 
be the $G$-orbit of the pair $(x(L_\emptyset),L_\emptyset)$.
There is an algebras isomorphism
$\widehat\Hb_{l|q^{-1}}{\buildrel\sim\over\to}\CC_G(X\times X)$
which maps $T_x$ to the characteristic function of $X_x$ (see [IM]). 
Put $P=X_1^{-1}\tilde T_1^{-1}\tilde T_2^{-1}\cdots\tilde T_{l-1}^{-1}$.

\vskip3mm

\noindent{\bf Lemma.} {\it The right action of $\widehat\Hb_l$ on
$\CC_G(Y\times X)$ is such that if $\ib\in A_l^n$, $x\in{\frak S}^\ib$,
and $s\in\widehat S_l$, then $(\1b_\ib)\,P=\otimes\1b_{(\ib)\pi}$ and
$$(\1b_{(\ib)x})\,\tilde T_s=\left\{\matrix 
v^{-1}\1b_{(\ib)x}\hfill&
\text{if}\quad xs\notin{\frak S}^\ib\quad(\text{then}\quad xs>x),\hfill\cr\cr
\1b_{(\ib)xs}\hfill&
\text{if}\quad xs>x\quad\text{and}\quad xs\in{\frak S}^\ib,\hfill\cr\cr
\1b_{(\ib)xs}+(v^{-1}-v)\1b_{(\ib)x}\hfill&
\text{if}\quad xs<x\quad(\text{then}\quad xs\in{\frak S}^\ib).\hfill\cr\cr
\endmatrix\right.$$
}

\vskip3mm

\noindent{\it Proof.} To simplify the notations fix $l=2$.
Fix $\ib\in A_l^n$ and $x\in{\frak S}^\ib$. Set $(i,j)=(\ib)x$. Then
$$xs_1\notin{\frak S}^\ib\iff
(\exists t\in S_\ib\quad\text{such\ that}\quad xs_1=tx)
\iff (i,j)s_1=(i,j).$$
Moreover,
$$xs_1>x\quad\text{and}\quad xs_1\in{\frak S}^\ib
\iff i\leq j\quad\text{and}\quad (i,j)s_1\neq (i,j).$$
The formulas in the proposition gives 
$$(\1b_{(i,j)})\,T_1=\left\{\matrix 
v^{-2}\1b_{(i,j)}\hfill&\text{if}\quad i=j,\hfill\cr\cr
v^{-1}\1b_{(j,i)}\hfill&\text{if}\quad i<j,\hfill\cr\cr
v^{-1}\1b_{(j,i)}+(v^{-2}-1)\1b_{(i,j)}\hfill&
\text{if}\quad i>j.\hfill
\endmatrix\right.$$
These are precisely the formulas in [VV, Section 5] taking into account
the different normalizations for the Hecke algebra and the factor 
$q^{-l(\omega_\ib)}$. The result follows from [VV, Proposition 6].\qed

\vskip3mm

\noindent{\bf 7.4.}
Fix $\ib\in A_l^n.$ Let $\Hb_\ib\subseteq\widehat\Hb_l$ 
be the parabolic subalgebra associated to ${\frak S}_\ib$. 
Set $e_\ib=\sum_{x\in{\frak S}_\ib}T_x$ and
$\pi_\ib=\sum_{x\in{\frak S}_\ib}v^{-2l(x)}.$
Thus $e_\ib^2=\pi_\ib e_\ib$ and $\bar e_\ib=v^{2l(\omega_\ib)}e_\ib$.
Set
$${\ts\Tb_{n,l}=\bigoplus_{\ib\in A_l^n}e_\ib\,\widehat\Hb_l.}$$
The affine $q$-Schur algebra $\widehat\Sb_{n,l}$, introduced in [G2],
is the endomorphism ring of the right $\widehat\Hb_l$-module
$\Tb_{n,l}.$ If $\jb\in A_l^n$ set
$$M_{\ib\jb}=\{\mb\in M^l\,|\,
Y_\mb\cap(G(L_\ib)\times G(L_\jb))\neq\emptyset\}.$$
A matrix $\mb\in\Mb_{\ib\jb}$ is identified with the class in
${\frak S}_\ib\setminus\widehat{\frak S}_l/{\frak S}_\jb$ 
of the elements $x$ such that 
$(L_{(\ib)x},L_\jb)\in Y_\mb$.
Set $T_\mb=\sum_{x\in\mb}T_x$. Let 
$\widehat\Hb_{\ib\jb}\subseteq\widehat\Hb_l$ be the $\AA$-linear span
of the elements $T_\mb$ with $\mb\in M_{\ib\jb}.$
The $\AA$-linear homomorphism 
$${\ts\bigoplus_{\ib,\jb\in A_l^n}\widehat\Hb_{\ib\jb}\to
\widehat\Sb_{n,l}}$$
which maps $T_\mb$, $\mb\in M_{\ib\jb}$, to the endomorphism such that
$e_\kb\mapsto\delta_{\kb,\jb}\,T_\mb\in e_\ib\widehat\Hb_l,$ is invertible. 
The product in the affine $q$-Schur algebra, denoted $\bullet$, is 
$$T_\mb\bullet T_{\nb}=\delta_{\kb,\jb}\,\pi_\jb^{-1}T_\mb T_\nb\qquad
\forall\mb\in M_{\ib\jb}\quad\forall\nb\in M_{\kb\lb}.$$
If $t\in\CC^\times$ let $\widehat\Sb_{n,l|t}$ and
$\Tb_{n,l|t}$ be the specializations of $\widehat\Sb_{n,l}$ and
$\Tb_{n,l}$ at $v=t$.

\vskip3mm

\noindent{\bf Proposition.} {\it $(a)$ The map 
$\Phi\,:\,\widehat\Sb_{n,l|q^{-1}}\to\CC_{G}(Y\times Y),$
$T_\mb\mapsto q^{y(\mb)}\1b_\mb$, 
is an isomorphism of algebras.\hfill\break
$(b)$ There is a unique isomorphism of
$\widehat\Sb_{n,l|q^{-1}}\times\widehat\Hb_{l|q^{-1}}$-modules, still denoted
$\Phi$, from $\Tb_{n,l|q^{-1}}$ to $\CC_G(Y\times X)$
such that $e_\ib\mapsto q^{l(\omega_\ib)}\1b_\ib$ for all $\ib\in A_l^n$.}

\vskip3mm

\noindent{\it Proof.} 
The map $\Phi$ is a linear isomorphism.
For all $x,y,z\in\widehat{\frak S}_l$ let $B_{xy}^z(v)\in\AA$ be such that
$${\ts T_xT_y=\sum_zB_{xy}^z(v)T_z.}$$
If $(L'',L)\in X_z$ then,
$$B_{xy}^z(q^{-1})=
\sharp\{L'\in X\,|\,(L'',L')\in X_x\quad\&\quad(L',L)\in X_y\}.$$
Fix $\mb\in M_{\ib\kb}$, $\nb\in M_{\ib\jb}$, and
$\ob\in M_{\jb\kb}$. Let $A_{\nb\ob}^{\mb}\in\NN$ be such that
$$\1b_{\nb}\star\1b_{\ob}=\sum_\mb q^{-y(\nb)-y(\ob)+y(\mb)}
A_{\nb\ob}^\mb\,\1b_\mb.$$
Then for any $z\in\mb$,
$${\ts A_{\nb\ob}^\mb=h_\jb^{-1}\sum_{x\in\nb\atop y\in\ob}
B_{xy}^z(q^{-1}),}$$
where $h_\jb=\pi_{\jb|v=q^{-1}}$ is the cardinal of the fiber of the projection 
$X\to G(L_\jb)$. Claim $(a)$ follows from the identity 
$${\ts T_{\nb}\bullet T_{\ob}=
\pi_\jb^{-1}\sum_z\sum_{x\in\nb\atop y\in\ob}B_{xy}^z(v)T_z=
\sum_\mb A_{\nb\ob}^\mb T_\mb\quad\mod (v-q^{-1}).}$$ 
Let $\Phi$ be the unique isomorphism of right 
$\widehat\Hb_l$-modules
$\Tb_{n,l|q^{-1}}{\buildrel\sim\over\to}\CC_G(Y\times X)$
such that $\Phi(e_\ib)=q^{l(\omega_\ib)}\1b_\ib$ for all $\ib\in A^n_l$.
Let us prove that $\Phi$ commutes to the action of $\widehat\Sb_{n,l}$.
We must prove that for all $\ib,\jb\in A^n_l$ and all 
$\mb\in M_{\ib\jb}$ then
$$q^{y(\mb)}\1b_\mb\star\1b_\jb=q^{-l(\omega_\jb)}h^{-1}_\ib\Phi(e_\ib T_\mb)=
q^{l(\omega_\ib)-l(\omega_\jb)}h^{-1}_\ib(\1b_\ib)\,T_\mb.$$
Put 
$$\1b_\mb\star\1b_\jb=\sum_{\kb\in\ZZ^l}
q^{-y(\mb)-l(\omega_\jb)+l(\omega_\kb)}A_{\mb\jb}^\kb\1b_\kb,\qquad
(\1b_\ib)T_\mb=\sum_{\kb\in\ZZ^l}q^{-l(\omega_\ib)+l(\omega_\kb)}
A_{\ib\mb}^\kb\1b_\kb.$$
Fix $z\in\widehat{\frak S}_l$ and $(L'',L)\in X_z$ whose projection in 
$Y\times X$ is in $X_\kb$. Then
$$A_{\mb\jb}^\kb=\sharp\{L'\in Y\,|\,(L'',L')\in Y_\mb,\quad
(L',L)\in X_\jb\}=
h_\jb^{-1}\sum_{y\in\mb\atop x\in{\frak S}_\jb}B_{yx}^z(q^{-1}),$$
$$A_{\ib\mb}^\kb=\sum_{y\in\mb}\sharp\{L'\in X\,|\,
(L'',L')\in X_\ib,\quad (L',L)\in X_y\}=
\sum_{y\in\mb\atop x\in{\frak S}_\ib}B_{xy}^z(q^{-1}).$$
Claim $(b)$ follows from the equality
$$\pi^{-1}_\ib\sum_{y\in\mb\atop x\in{\frak S}_\ib}T_xT_y=
\pi^{-1}_\jb\sum_{y\in\mb\atop x\in{\frak S}_\jb}T_yT_x.$$
\qed

\vskip3mm

\noindent{\bf 7.5.}
The set $M^+=\{\mb\in M\,|\,i>j\Rightarrow m_{ij}=0\}$
parametrizes the isomorphism classes of nilpotent
representations of the quiver $\Gamma_n$ : $O_\mb$ is the class of 
$\bigoplus_{i=1}^n\bigoplus_{j\geq i}\bar\FF[i,j]^{m_{ij}}$.
Let $\ \bar{}\ $ be the unique semilinear involution on 
$\Ub^-_n$ fixing the elements $\bb_{O_\mb}$.

\vskip3mm

\noindent{\bf Proposition.} {\it The involution $\ \bar{}\ $ on $\Ub^-_n$
is a ring homomorphism and $\bar\fb_\alpha=\fb_\alpha$ for all
$\alpha\in\NN^{(\ZZ/n\ZZ)}$.}

\vskip3mm

\noindent{\it Proof.} The second claim is obvious since $\fb_\alpha$
is the characteristic function of a single point. We now prove
the first claim. For any algebraic variety $X$ over $\bar\FF$ let
$\Dc(X)$ be the bounded derived category of complexes of $\QQ_l$-sheaves
on $X$ (see [L2], [L3]). If $G$ is a connected algebraic group
acting on $X$ let $\Dc_G^{ss}(X)$ be the full subcategory whose objects
are sums of shifted simple $G$-equivariant objects in $\Dc(X)$. 
Lusztig has defined in [L2, Section 3.1] a convolution product 
$$*\,:\,\Dc_{G_U}^{ss}(E_U)\times\Dc_{G_W}^{ss}(E_W)\to\Dc_{G_V}^{ss}(E_V)$$
such that $\Fc*\Gc=(p_3)_!\Hc$ where $\Hc$ satisfies 
$p_2^*\Hc\simeq p_1^*(\Fc\otimes\Gc)$. Let $D$ be the Verdier duality.
Since $p_1$ and $p_2$ are smooth with connected fibers
and since $p_3$ is proper we get $D(\Fc*\Gc)=(D\Fc)*(D\Gc)[2d_1-2d_2]$ 
where $d_1$ and $d_2$ are the dimensions of the fibers of $p_1$ and $p_2$.
Let $\alpha,\beta$ be the dimension of $U,W$.
We know that $d_2=\sum_\io\alpha_\io^2+\sum_\io\beta_\io^2$ and 
$d_1=d_2+m(\beta,\alpha)$. Thus $D(\Fc*\Gc)=(D\Fc)*(D\Gc)[2m(\beta,\alpha)]$.
Finally observe that the elements $\bb_{O_\mb}$ are the Frobenius traces
of the simple perverse sheaves on the $E_V$ since the varieties
$\bar O_\mb$ are pure (see [L1, Corollary 11.6]).
\qed

\vskip3mm

\noindent{\bf 7.6.}
If $L',L\in Y$ are such that $L'\subseteq L$ then $L/L'$
may be viewed as a nilpotent representation of $\Gamma_n$ of
dimension $\alpha$ where $\alpha_\io=\dim_{\FF}(L_i/L'_i)$ 
(see [L2, Section 11], [GV]). Then, set 
$$a(L',L)=\sum_{i=1}^n
\dim_{\FF}(L_i/L'_i)(\dim_{\FF}(L'_{i+1}/L'_i)-\dim_{\FF}(L_i/L'_i)).$$
Let $\Theta\,:\,\Ub^-_n\to\widehat\Sb_{n,l}$
be the $\AA$-linear map such that 
$$\Phi\circ\Theta(f)(L',L)=q^{-a(L',L)}f(L/L')\quad\text{if}\quad
L'\subseteq L,\quad 0\quad\text{else}.$$
If $\ib\in A^n_l$ and $\mb\in M^+$ let $\mb^\ib\in\cup_\jb M_{\ib\jb}$
be the matrix with the $(i,j)$-th entry equal to
$${\ts\delta_{ij}(\sharp\ib^{-1}(j+1)-\sum_{k\leq j}m_{kj})+
(1-\delta_{ij})m_{i+1,j}.}$$
Let $\phi$ be the semilinear involution on $\widehat\Sb_{n,l}$ such that
$\phi(u)=v^{-2l(\omega_\jb)}\bar u$ for all $u\in\widehat\Hb_{\ib\jb}.$

\vskip3mm

\noindent{\bf Proposition.} {\it 
The map $\Theta\,:\,\Ub^-_n\to\widehat\Sb_{n,l}$
is an algebra homomorphism. Moreover if
$u\in\Ub^-_n$ and $\mb\in M^+$ we have
$$\phi\circ\Theta(u)=\Theta(\bar u)
\quad\text{and}\quad
\Phi\circ\Theta(\fb_{O_\mb})=\sum_{\ib\,|\,\mb^\ib\in M}\1b_{\mb^\ib}.$$
}

\vskip3mm

\noindent{\it Proof.} The first claim is immediate from the formula
$$a(L'',L)-a(L',L)-a(L'',L') =-m(L/L',L'/L'')$$
and from the definition of the product in 
$\Ub^-_n$ and $\CC_G(Y\times Y)$.
We know that $\bar\fb_\alpha=\fb_\alpha$ for all $\alpha\in\NN^{(\ZZ/n\ZZ)}$.
Similarly, $\phi\circ\Theta(\fb_\alpha)=\Theta(\fb_\alpha)$ 
since for any flag $L'$ the $L$'s such that $L'\subseteq L$ and
$\fb_\alpha(L/L')\neq 0$ are the rational points of a smooth
variety. Hence, the second claim results from
the first claim, Proposition 3.5, and the fact that $\phi$ is a ring
homomorphism. Now let us first prove that 
$$(L',L)\in Y_{\mb^\ib}\quad\iff\quad
\biggl(L/L'\in O_\mb\quad\text{and}\quad
\dim_{\FF}(L'_i/L'_{i-1})=\sharp\ib^{-1}(i)\biggr).\leqno{(c)}$$
By definition $(L',L)\in Y_{\mb^\ib}$ if and only if
$${\ts\delta_{ij}(\sharp\ib^{-1}(j+1)-\sum_{k\leq j}m_{kj})+
(1-\delta_{ij})m_{i+1,j}
=\dim_{\FF}\biggl({L_{i+1}\cap L'_{j+1}\over
(L_i\cap L'_{j+1})+ (L_{i+1}\cap L'_j)}\biggr).}$$
Thus it suffices to prove that if
$\dim_{\FF}(L'_{i+1}/L'_i)=\sharp\ib^{-1}(i+1)$ and $L'\subseteq L$ then 

\vskip2mm

\noindent $(d)\quad\dim_{\FF}\biggl({L_i\cap L'_{j+1}\over
(L_{i-1}\cap L'_{j+1})+ (L_i\cap L'_j)}\biggr)$ is the multiplicity 
of $\FF[i,j]$ in $L/L'$ for all $i\leq j$,

\vskip2mm

\noindent $(e)\quad$if $x_\io\,:\,L_i/L'_i\to L_{i+1}/L'_{i+1}$ is the map
induced by the inclusion $L_i\subseteq L_{i+1}$, then
$$\sharp\ib^{-1}(i+1)-\dim_{\FF}\Ker(x_\io)=\dim(L'_{i+1}/(L_i\cap L'_{i+1})).$$

\noindent Claim $(e)$ is immediate since
$$\Ker(x_\io)=(L_i\cap L'_{i+1})/L'_i,\qquad
L'_{i+1}/(L_i\cap L'_{i+1})\simeq 
{L'_{i+1}/L'_i\over (L_i\cap L'_{i+1})/L'_i},$$
and since $\sharp\ib^{-1}(i+1)=\dim(L'_{i+1}/L'_i)$. 
Part $(d)$ is due to the fact that
$\FF[i,j]$ is a direct summand of $L/L'$ 
if and only if there is a vector $w\in L'_{j+1}\setminus L'_j$ such that
$w\in L_i\setminus L_{i-1}$. The second claim follows from $(c)$ and 
the formula
$$(\mb\in M^+\quad\text{and}\quad\mb^\ib\in M)\Rightarrow
\dim O_\mb+a(\mb^\ib)=y(\mb^\ib),$$
which is left to the reader.\qed

\vskip3mm

\noindent{\bf Remark.} 
For any $\mb\in M^l$ set
$$\cb_\mb=\sum_{i,\nb}v^{-i+y(\mb)}
\dim\Hc^{i}_{Y_{\nb,L}}(IC_{Y_{\mb,L}})\,T_{\nb},$$
where $L$ is such that $Y_{\mb,L}\neq\emptyset$.
The elements $\cb_\mb$ form a $\AA$-basis of $\widehat\Sb_{n,l}$
and $\phi(\cb_{\mb})=\cb_\mb$. 
Proposition 7.6 implies that for any $\mb\in M^+$ we have
$$\Theta(\bb_{O_\mb})=\sum_{\ib\,|\,\mb^\ib\in M}\cb_{\mb^\ib}.$$
We will not use this.

\vskip5mm

\centerline{\bf 8. The tensor representation of $\tilde\Ub^-_n$.}

\vskip3mm

\noindent{\bf 8.1.}
Let $\AA^{(\ZZ)}$ be the $\AA$-linear span of vectors $\x_i$, $i\in\ZZ$.
Let $\eb_{ij}\in M$ be the matrix with $1$ 
at the spot $(k,l)$ if $(k,l)\in (i,j)+\ZZ\,(n,n)$ and $0$ elsewhere. 

\vskip3mm

\noindent{\bf Lemma.} {\it $\tilde\Ub^-_n$ acts on $\AA^{(\ZZ)}$ in such 
a way that for all $\mb\in M$ and all $\alpha\in\NN^{\ZZ/n\ZZ}$,
$$\fb_{O_\mb}(\x_i)=\sum_{j\geq i}\delta_{\mb,\eb_{ij}}\,\x_{j+1}
\qquad\text{and}\qquad
\kb_\alpha(\x_i)=v^{-n(\alpha,\epsilon_\io)}\x_i.$$
}

\noindent{\it Proof.} It is the action obtained by taking
$l=1$ in the geometric construction of Section 7 via the isomorphism
$\Tb_{n,1}{\buildrel\sim\over\to}\AA^{(\ZZ)}$, $\1b_i\mapsto\x_i$
(if $l=1$ then $X$ and $Y$ are zero dimensional).
\qed

\vskip3mm

\noindent{\bf 8.2.}
Put $\bigotimes^l=(\AA^{(\ZZ)})^{\otimes l}$.
For any sequence $\ib=(i_1,i_2,...,i_l)\in\ZZ^l$ set
$\otimes\x_\ib=\x_{i_1}\otimes\x_{i_2}\otimes\cdots\otimes\x_{i_l}.$
On one hand $\bigotimes^l$ is a left $\tilde\Ub_n^-$-module
via the coproduct $\Delta$. On the other hand
$\widehat\Hb_l$ acts on $\bigotimes^l$ as follows for all
$k=1,2,...,l-1$ and $j=1,2,...,l$ :

$$(\otimes\x_\ib)T_k=\left\{\matrix 
v^{-2}\otimes\x_\ib\hfill&\text{if}\quad i_k=i_{k+1}\hfill\cr\cr
v^{-1}\otimes\x_{(\ib)s_k}\hfill&\text{if}\quad -n<i_k<i_{k+1}\leq 0
\hfill\cr\cr
v^{-1}\otimes\x_{(\ib)s_k}+(v^{-2}-1)\otimes\x_\ib\hfill&
\text{if}\quad -n<i_{k+1}<i_k\leq 0,\hfill
\endmatrix\right.\leqno(a)$$

\noindent{$(b)\qquad(\otimes\x_\ib)X_j^{-1}=\otimes\x_{(\ib)\epsilon_j}$.

\vskip3mm

\noindent{\bf Lemma.} {\it The representations of $\tilde\Ub^-_n$ and
$\widehat\Hb_l$ on $\bigotimes^l$ commute.}

\vskip3mm

\noindent{\it Proof.} Since the coproduct is coassociative 
(see [G1, Theorem 1(ii)]), we are reduced to the case $l=2$. By definition
$$\Delta(\fb_\alpha)=\sum_{\alpha=\beta+\gamma}v^{n(\gamma,\beta)}
\fb_{\beta}\kb_{\gamma}\otimes\fb_{\gamma}.\leqno(c)$$
Thus $\fb_\alpha$ acts on $\bigotimes^2$ as
$$\left\{\matrix
\fb_\io\otimes 1+\kb_\io\otimes\fb_\io\hfill&
\text{if}\quad\alpha=\epsilon_\io\hfill\cr\cr
\fb_\io\otimes\fb_\jo+\fb_\jo\otimes\fb_\io\hfill&
\text{if}\quad\alpha=\epsilon_\io+\epsilon_\jo\quad\text{and}\quad\io\neq\jo
\hfill\cr\cr
\fb_\io\otimes\fb_\io\hfill&
\text{if}\quad\alpha=2\epsilon_\io\hfill\cr\cr
0\hfill&\text{else}.\hfill
\endmatrix\right.$$
The commutation results from a direct computation.
\qed

\vskip3mm

\noindent{\bf 8.3. Lemma.} {\it The map 
$\otimes\x_\ib\mapsto v^{l(\omega_\ib)}e_\ib$, for all $\ib\in A_l^n$, 
extends uniquely to an isomorphism of 
$\Ub^-_{n}\times\widehat\Hb_l$-bimodules 
$\bigotimes^l{\buildrel\sim\over\to}\Tb_{n,l}$.
}

\vskip3mm

\noindent{\it Proof.}
The map above extends uniquely to an isomorphism of $\widehat\Hb_l$-modules.
Let us prove that this isomorphism commutes to $\Ub^-_{n}$.
For any $\ib\in\ZZ^l$ Lemma $8.1$ and formula $(c)$ give
$${\ts\fb_\alpha(\otimes\x_\ib)=
\sum_\nb v^{c(\ib,\ib+\nb)}\otimes\x_{\ib+\nb},}$$
where $\nb=(n_1,n_2,...,n_l)\in\{0,1\}^l$ describes the set of all 
sequences such that $\alpha=\sum_{s=1}^ln_s\epsilon_{\io_s}$ and
$$c(\ib,\ib+\nb)=-\sum_{1\leq s<t\leq l}n_t(1-n_s)
n(\epsilon_{\io_t},\epsilon_{\io_s}).$$
By Proposition 7.4, after the specialisation
$v=q^{-1}$ we have 
$q^{-l(\omega_\ib)}\fb_\alpha(e_\ib)=(\Phi\circ\Theta)(\fb_\alpha)\star\1b_\ib$.
The R.H.S. is the convolution product of $\1b_\ib$ and a function supported
by the set of all pairs $(L',L)$ such that for all $i$ we have
$$L'_i\subseteq L_i\subseteq L'_{i+1},\qquad
\dim_\FF(L_i/L'_i)=\alpha_i,
\quad\text{and}\quad \dim_\FF(L_i/L_{i-1})=\sharp\ib^{-1}(i).\leqno(d)$$
By definition of the convolution product 
$(\Phi\circ\Theta)(\fb_\alpha)\star\1b_\ib$ 
(simply denoted by $\fb_\alpha(\1b_\ib)$) is a linear 
combination of the $\1b_\jb$'s such that it exists $L\in Y$ such that 
$(L,L_\emptyset)\in X_\ib$ and $(L_\jb,L)$ satisfies $(d)$.
Suppose first that $\ib\in A^n_l$. Then
$X_\ib\cap(Y\times\{L_\emptyset\})=\{(L_\ib,L_\emptyset)\}.$
Thus $\fb_\alpha(\1b_\ib)$ is a linear combination of the $\1b_\jb$'s such that 
$$\ib\leq\jb\leq\ib+1\quad\text{and}\quad
\alpha_i=\sharp(\ib^{-1}(i)\cap\jb^{-1}(i+1)).\leqno(e)$$
More precisely we get
$$\fb_\alpha(\1b_\ib)=\sum_\nb 
q^{-a(\ib+\nb,\ib)-l(\omega_\ib)+l(\omega_{\ib+\nb})}\1b_{\ib+\nb},$$
where the $\nb$'s are as above and
$a(\ib+\nb,\ib)=\sum_i\alpha_i(\sharp\ib^{-1}(i+1)-\alpha_{i+1})$. 
A simple computation using $(e)$ gives
$$a(\ib+\nb,\ib)=\sum_{s,t=1}^ln_t(1-n_s)\delta_{\io_t+1,\io_s}.$$
Moreover for any sequence $\jb$ we have 
$l(\omega_\jb)=\dim(P_\jb/(B\cap P_\jb))$ 
where $P_\jb,B\subset G$ are the isotropy subgroup of $L_\jb$ and $L_\emptyset$.
Thus we obtain
$$l(\omega_{\ib+\nb})-l(\omega_\ib)=
\sum_{1\leq t<s\leq l}n_t(1-n_s)\delta_{\io_s,\io_t+1}(1-\delta_{i_t,0})+$$
$$+\sum_{1\leq t<s\leq l}n_s(1-n_t)\delta_{i_t,1-n}\delta_{i_s,0}-
\sum_{1\leq t<s\leq l}n_s(1-n_t)\delta_{\io_t,\io_s}.$$
To conclude it suffices to compute the image of $\otimes\x_{\ib+\nb}$
in $\CC_G(Y\times X)$. Using the identity 
$$X_k^{-1}=\tilde T_{k-1}^{-1}\cdots\tilde T_1^{-1}P\tilde T_{l-1}\cdots
\tilde T_k$$
and Lemmas 7.3 and 8.2, we get that $\otimes\x_{\ib+\nb}$ is mapped to
$$q^{-d(\ib,\ib+\nb)}\1b_{\ib+\nb},\quad
d(\ib,\ib+\nb)=\sharp\{1\leq s,t\leq l\,|\,i_s=0,\,i_t=1-n,\,n_s=1,\,n_t=0\}.$$
Then, the equality results from an easy computation.
The general case (i.e. $\ib\notin A^n_l$) follows since the isomorphism we
consider commutes to $\widehat\Hb_l$.
\qed

\vskip3mm

\noindent{\bf 8.4.}
Let $\psi$ be the semilinear involution on $\Tb_{n,l}\simeq\bigotimes^l$ 
such that $\psi(e_\ib t)=\bar e_\ib\bar t$ for all $t\in\widehat\Hb_l.$
Proposition 7.6 and the definition of the involutions $\phi$ and $\psi$ 
imply the following lemma (see Subsection 7.6).

\vskip3mm

\noindent{\bf Lemma.} {\it 
For all $u\in\Ub^-_n$ and all $t\in\bigotimes^l$ we have
$\psi(u\, t)=\bar u\psi(t)$.}
\qed

\vskip5mm

\centerline{\bf 9. The action of $\Ub^-_n$ on wedges.}

\vskip3mm

\noindent{\bf 9.1.}
Set $\Omega^l=\sum_i\Im(T_i+1)\subset\bigotimes^l$.
We have
$${\ts\bigotimes^l/\Omega^l}\simeq
\bigoplus_{\ib\in A_l^n} e_\ib\,\widehat\Hb_l\,e^-,$$
where $e^-=\sum_{x\in{\frak S}_l}(-v)^{l(x)}\tilde T_x$. 
For any $\ib\in\ZZ^l$ let $\wedge\x_\ib$ be the image of $\otimes\x_\ib$
in $\bigotimes^l/\Omega^l$. Set
$$P_l^{++}=\{\ib\in\ZZ^l\,|\,i_1>i_2>...>i_l\}.$$
The $\wedge\x_\ib$'s such that $\ib\in P_l^{++}$ 
form a basis of $\bigotimes^l/\Omega^l$ (see [KMS, Proposition 1.3]).
For any $\lambda\in\Pi_l$ set
$|\lambda\rangle=\wedge\x_\ib$ if $\ib=\lambda+\rho$, where
$\rho$ is as in Section 1.1. 
Let $\bigwedge^l\subset\bigotimes^l/\Omega^l$ be the linear
span of the vectors $|\lambda\rangle.$

\vskip3mm

\noindent{\bf 9.2.}
The representation of $\Ub^-_n$ on $\bigotimes^l$
descends to $\bigotimes^l/\Omega^l$ (use Lemma 8.2).
For all $\lambda\in\Pi_l$ set 
$\bb_\lambda=\bb_{O_\lambda}|\emptyset\rangle$ and put
$\Bb_l=\{\bb_\lambda\,|\,\lambda\in\Pi_l\}.$
Let consider the involution $\psi$ on $\bigotimes^l/\Omega^l$
such that
$$\psi(e_\ib\,h\,e^-)=v^{2l(\omega)}\bar e_\ib\,\bar h\,e^-
\qquad\forall h\in\widehat\Hb_l.$$

\vskip3mm

\noindent{\bf Proposition.} {\it  $\Bb_l$ is a basis of $\bigwedge^l$ 
whose elements are fixed by $\psi$.}

\vskip3mm

\noindent{\it Proof.} 
Lemma $8.1$, the definition of $\Delta$, and the normal ordering rule 
[KMS, $(43)$ and $(45)$] imply that for any $\lambda\in\Pi_l$ 
and any orbit $O\subset\bar O_\lambda\setminus O_\lambda$ we have
$${\ts
\fb_{O_\lambda}(|\emptyset\rangle)\in v^\ZZ|\lambda\rangle+
\bigoplus_{\mu<\lambda}\AA|\mu\rangle
\quad\text{and}\quad
\fb_O(|\emptyset\rangle)\in\bigoplus_{\mu<\lambda}\AA|\mu\rangle.}$$
Thus, $\Bb_l$ is a basis.
Now Lemma 8.4 implies that the action of $\bb_{O_\mb}$ on 
$\bigwedge^l$ commutes to $\psi$.
Since $\psi(|\emptyset\rangle)=|\emptyset\rangle$ we get
$\psi(\bb_\lambda)=\bb_\lambda$ for all $\lambda$.
\qed

\vskip3mm

\noindent{\bf 9.3.}
Let $P^+_l\subset\ZZ^l$ be the subset of integral dominant weights.
Let ${\frak S}^{\ib,l}$ be the set of minimal length 
representatives of the cosets in 
${\frak S}_\ib\setminus\widehat{\frak S}_l/{\frak S}_l$. 
Thus ${\frak S}^{\ib,l}={\frak S}^\ib\cap{\frak S}^l$
where ${\frak S}^l$ is the set of minimal length 
representatives of the cosets in $\widehat{\frak S}_l/{\frak S}_l$. 
For any $x\in\widehat{\frak S}_l$ let $\check x$ be the smallest element
in the double coset ${\frak S}_\ib x{\frak S}_l$. Set
$${\frak S}(\ib,l)=\{x\in{\frak S}^{\ib,l}\,|\,S_\ib x\cap x S_l=\emptyset\}.$$

\vskip3mm

\noindent{\bf Lemma.} {\it Fix $\ib\in A_l^n$. Then,\hfill\break
$(a)\quad e_\ib\tilde T_x e^-\neq 0\Rightarrow
x\in{\frak S}_\ib{\frak S}(\ib,l){\frak S}_l,$\hfill\break
$(b)\quad(x\in{\frak S}^\ib\quad\&\quad(\ib)x\in P^{++}_l)\Rightarrow
e_\ib\tilde T_xe^-=v^{-l(\omega_\ib)}\wedge\x_{(\ib)x},$\hfill\break
$(c)\quad(\ib)x\in P^{++}_l\iff x\in{\frak S}_\ib{\frak S}(\ib,l)\omega$.
}

\vskip3mm

\noindent{\it Proof.}
Suppose that $x\in{\frak S}^{\ib,l}$, $s_i\in S_\ib$, and $s\in S_l$,
are such that $s_ix=xs$. Then
$$v^{-1}e_\ib\tilde T_xe^-=e_\ib\tilde T_{s_i}\tilde T_xe^-=
e_\ib\tilde T_x\tilde T_se^-=-v e_\ib\tilde T_xe^-.$$
Thus $e_\ib\tilde T_xe^-=0$. Any $x\in\widehat{\frak S}_l$  
decomposes in $x=x_\ib \check x x_l$ where $x_\ib\in{\frak S}_\ib$, 
$\check x\in{\frak S}^{\ib,l},$ $x_l\in{\frak S}_l$, and $l(x)=l(x_\ib)+l(\check x)+l(x_l)$. In particular 
$$e_\ib\tilde T_xe^-=v^{-l(x_\ib)}(-v)^{l(x_l)}
e_\ib\tilde T_{\check x}e^-.$$
Claim $(a)$ follows. Let us prove claim $(b)$.
Recall that if $\lambda$ is dominant then $\tilde T_\lambda^{-1}$ is mapped to 
$X^\lambda=X_1^{\lambda_1}X_2^{\lambda_2}\cdots X_l^{\lambda_l}$ by the 
Bernstein isomorphism. 
Then $(8.2.b)$ implies that
$$(\otimes\x_\ib)\tilde T_\lambda=\otimes\x_{(\ib)\lambda}\qquad
\forall\lambda\in P_l^+\quad\forall\ib\in\ZZ^l,$$
Moreover $(8.2.a)$ implies that
$$(\otimes\x_\ib)\tilde T_x=\otimes\x_{(\ib)x}\qquad
\forall x\in{\frak S}_l\cap{\frak S}^\ib\quad\forall\ib\in A_l^n.$$
Fix $x\in\widehat{\frak S}_l$. Then $x$ decomposes uniquely as 
$x=y\lambda$ where $y\in{\frak S}_l$ and $\lambda\in\ZZ^l$.
If $(\ib)x=(\ib)y+n\lambda\in P^{++}_l$ then $\lambda\in P_l^+$. 
Since $s\lambda>\lambda$ for all $s\in S_l$ if $\lambda$ is dominant, we get 
$\tilde T_x=\tilde T_y\tilde T_\lambda$.
Suppose moreover that $x\in{\frak S}^\ib$. Then for any $s\in S_\ib$
we have $sy\lambda>y\lambda$. Since $l(z\lambda)=l(z)+l(\lambda)$ for
any $z\in{\frak S}_l$ ($\lambda$ is dominant), we obtain that 
$y\in{\frak S}_l\cap{\frak S}^\ib$. Hence, Section 8.2 implies that
$$e_\ib\tilde T_xe^-=e_\ib\tilde T_y\tilde T_\lambda e^-=
v^{-l(\omega_\ib)}\wedge\x_{(\ib)x}.$$
Finally, claim $(c)$ follows from
$$\matrix
(x\in{\frak S}^\ib\quad\&\quad(\ib)x\in P^{++}_l)&\iff
(s_ix>x>xs\quad\forall s\in S_l\quad\forall s_i\in S_\ib)\hfill\cr\cr
&\iff x\in{\frak S}(\ib,l)\omega.\hfill\endmatrix$$
\qed

\vskip3mm

\noindent{\bf Proposition.} {\it If $\ib\in A_l^n$ and $(\ib)x\in P^{++}_l$ then
$\psi(\wedge\x_{(\ib)x})=(-1)^{l(\omega)}v^{l(\omega^\ib)}
\wedge\x_{(\ib)x\omega}.$}

\vskip3mm

\noindent{\it Proof.} 
First recall that if 
$\lambda\in P^+_l$ and $\lambda^*=-\omega(\lambda)$, then 
$\tilde T_{\lambda^*}=\tilde T_\omega^{-1}\tilde T_{-\lambda}\tilde T_\omega$
(indeed, since $\lambda,\lambda^*$ are dominant weights we have 
$\tilde T_\omega\tilde T_{\lambda^*}=\tilde T_{-\lambda\omega}=
\tilde T_{-\lambda}\tilde T_\omega.$ In particular, 
$\tilde T_{-\lambda}=\tilde T_\omega\tilde T_{\lambda^*}\tilde T_\omega^{-1}$).
Fix $\ib\in A^n_l$ and $x\in{\frak S}^\ib$ such that $(\ib)x\in P^{++}_l$. 
As above fix $x=y\lambda$ with 
$y\in{\frak S}_l\cap{\frak S}^\ib$ and $\lambda\in P^+_l$. 
Using Lemma $9.3.b$ we get
$$\psi(\wedge\x_{(\ib)x})=\psi((\otimes\x_\ib)\tilde T_y\tilde T_\lambda e^-)=
v^{-l(\omega_\ib)}\psi(e_\ib\tilde T_y\tilde T_\lambda e^-)=
v^{l(\omega_\ib)+2l(\omega)} e_{\ib}\tilde T^{-1}_{y^{-1}}
\tilde T^{-1}_{-\lambda}e^-=$$
$$=v^{l(\omega_\ib)+2l(\omega)}e_{\ib}\tilde T^{-1}_{y^{-1}}
\tilde T_\omega\tilde T^{-1}_{\lambda^*}\tilde T_\omega^{-1}e^-=
(-1)^{l(\omega)}v^{2l(\omega)-l(\omega^\ib)}
e_{\ib}\tilde T_{y\omega}\tilde T^{-1}_{\lambda^*}e^-.$$
For all $s\in S_\ib$, $sy>y$ implies that $sy\omega<y\omega$. Thus
$$\psi(\wedge\x_{(\ib)x})=
(-1)^{l(\omega)}v^{l(\omega^\ib)}(\otimes\x_{(\ib)y\omega})
\tilde T^{-1}_{\lambda^*}e^-=
(-1)^{l(\omega)}v^{l(\omega^\ib)}\wedge\x_{(\ib)x\omega}.$$
\qed

\vskip3mm

\noindent{\bf 9.4.}
Fix $x\in{\frak S}^l$. The element
$$D_{x\omega}=\sum_{y\in{\frak S}^l\atop y\leq x}
(-v)^{l(y)-l(x)}\bar P_{y\omega,x\omega}\tilde T_{y\omega}e^-$$
is fixed by the involution on $\widehat\Hb_{l}e^-$ such that 
$h e^-\mapsto v^{2l(\omega)}\bar h e^-$ and 
$${\ts D_{x\omega}-\tilde T_{x\omega} e^-
\in\bigoplus_{y\in{\frak S}_l}v^{-1}\bar\SS\,\tilde T_{y\omega} e^-}$$
(see [D1]).
If $\ib\in A^n_l$ and $x\in{\frak S}(\ib,l)$ 
set $\bb^-_{(\ib)x\omega}=v^{l(\omega_\ib)} e_\ib D_{x\omega}.$
Lemma $9.3$ gives
$$\bb^-_{(\ib)x\omega}=\sum_{y\in{\frak S}^l\atop y\leq x}
(-v)^{l(y)-l(x)}v^{-l(y_\ib)}\bar P_{y\omega,x\omega}
\,\wedge\x_{(\ib)y\omega}.$$
Similarly fix $\ib\in A^n_l$ and $x\in{\frak S}^\ib$. Then
$$C'_{\omega_\ib x}=v^{l(x)+l(\omega_\ib)}
\sum_{y\in{\frak S}^\ib\atop y\leq x}\sum_{z\in{\frak S}_\ib}
P_{zy,\omega_\ib x}\,T_{zy}.$$
If $y\in{\frak S}^\ib$ and $y\leq x$ then 
$P_{zy,\omega_\ib x}=P_{\omega_\ib y,\omega_\ib x}$
for all $z\in{\frak S}_\ib$ (see [D1, page 491]).
Thus
$$C'_{\omega_\ib x}=v^{l(x)+l(\omega_\ib)}e_\ib
\sum_{y\in{\frak S}^\ib\atop y\leq x}P_{\omega_\ib y,\omega_\ib x}\,T_y.$$
If $x\in{\frak S}(\ib,l)$ set 
$\bb^+_{(\ib)x\omega}=(-v)^{l(\omega)}C'_{\omega_\ib x}e^-$.
Then Lemma $2.2$ gives
$$\matrix
\bb^+_{(\ib)x\omega}&=
\sum_{(y,z)}v^{l(x)-l(yz)}(-v)^{l(z)}
P_{\omega_\ib yz,\omega_\ib x}\,\wedge\x_{(\ib)y\omega}\hfill\cr\cr
&=\sum_{y\in{\frak S}(\ib,l)\atop y\leq x}v^{l(x)-l(y)}
Q_{\omega_\ib y,\omega_\ib x}\,\wedge\x_{(\ib)y\omega},\hfill
\endmatrix$$
where the first sum is over all couples
$(y,z)\in{\frak S}(\ib,l)\times{\frak S}_l$
such that $yz\leq x$ and $Q_{\omega_\ib y,\omega_\ib x}=
\sum_z(-1)^{l(z)}P_{\omega_\ib yz,\omega_\ib x}$ is a parabolic
Kazhdan-Lusztig polynomial.
Observe that $\bb^\pm_{(\ib)x\omega}$ is completely characterized by the 
following properties :
$$\psi(\bb^\pm_{(\ib)x\omega})=\bb^\pm_{(\ib)x\omega},\qquad
\bb^-_{(\ib)x\omega}-\wedge\x_{(\ib)x\omega}\in
\bigoplus_{y\in{\frak S}(\ib,l)\atop y<x}
v^{-1}\bar\SS\wedge\x_{(\ib)y\omega},$$
$$\text{and}\quad
\bb^+_{(\ib)x\omega}-\wedge\x_{(\ib)x\omega}\in
\bigoplus_{y\in{\frak S}(\ib,l)\atop y<x}
v\,\SS\wedge\x_{(\ib)y\omega}.$$
In particular  $\{\bb_\ib^-\,|\,\ib\in P^{++}_l\}$ and $\{\bb_\ib^+\, |\, \ib
\in P^{++}_l\}$ are bases of
$\bigotimes^l/\Omega^l$.  For all $\lambda\in\Pi_l$ set $\bb_\lambda^\pm=
\bb_\ib^\pm$ if $\ib=\lambda+\rho.$ Put $\Bb_l^\pm=\{\bb_\lambda^\pm\,
|\, \lambda\in\Pi_l\}.$

\vskip3mm

\noindent{\bf Remark.} If $\ib\in A^n_l$ and
$x,y\in{\frak S}^\ib$ are such that if $(\ib)x, (\ib)y\in P^{++}_l$ then
$$y\leq x\Rightarrow (\ib)x-(\ib)y\text{\ is\ a\ positive\ root}.$$

\vskip3mm

\noindent{\bf 9.5.}
The space $\bigwedge^l$ is endowed with four bases :
$\Bb^\pm_l=\{\bb^\pm_\lambda\,|\,\lambda\in\Pi_l\}$, 
$\Bb_l=\{\bb_\lambda\,|\,\lambda\in\Pi_l\}$, and
$\{|\lambda\rangle\,|\,\lambda\in\Pi_l\}.$
Moreover, $\Bb^\pm_l$ are characterized by
$$\psi(\bb^\pm_\lambda)=\bb^\pm_\lambda,\qquad
\bb^-_\lambda-|\lambda\rangle\in\bigoplus_{\mu<\lambda}
v^{-1}\bar\SS |\mu\rangle
\quad\text{and}\quad
\bb^+_\lambda-|\lambda\rangle\in\bigoplus_{\mu<\lambda}
v\,\SS |\mu\rangle$$
(in particular $\psi(\bigwedge^l)=\bigwedge^l).$  
Recall that if  $x\in\widehat{\frak S}_l$ and $\lambda\in\ZZ^l,$ then 
$\lambda\cdot x=(\lambda+\rho)x-\rho$ (see 1.1).
Section 9.4 implies the following theorem.

\vskip3mm

\noindent{\bf Theorem.} {\it 
$(a)\quad$If $\lambda\in\Pi_l$ and $x$ is minimal such that
$\ib=\lambda\cdot x^{-1}+\rho\in A^n_l$, then
$$\bb^-_{\lambda}=\sum_y(-v)^{l(y)-l(x)}v^{-l(y_\ib)}\bar P_{yx}
|{\lambda\cdot x^{-1}y\rangle},$$
where the sum is over all the $y$
such that $y\leq x$ and $\lambda\cdot x^{-1}y\in\Pi_l$.\hfill\break
$(b)\quad$For all $\lambda\in\Pi_l$ the coordinates of $\bb_{\lambda}^+$ 
in the wedges are some parabolic Kazhdan-Lusztig polynomials
(w.r.t. the parabolic subgroup ${\frak S}_l\subset\widehat{\frak S}_l$).
}
\qed

\vskip3mm

\noindent
Now, suppose that $l\leq n$. Let consider $\ib,\jb\in A^n_l$, 
$\nb\in M^+$ and $\mb=\nb^\jb\in M_{\jb\ib}\cap M^+$. 
If $x\in{\frak S}^\ib$ is such that $(\ib)x\in P_l^{++}$
then Section $7$ gives
$$\fb_{O_\nb}(\wedge\x_{(\ib)x})=v^{l(\omega_\ib)+y(\mb)}T_\mb\tilde T_x e^-
\in e_\jb\widehat\Hb_l e^-.$$
In particular if $\ib=(\rho)\omega\in A^n_l$ and $x=\omega$ we get
$$\fb_{O_\nb}(\wedge\x_\emptyset)=v^{l(\omega_\jb)}
e_\jb\tilde T_m\tilde T_\omega e^-,$$
where $m\in\mb$ is the smallest element.
Let suppose that $\fb_{O_\nb}(\wedge\x_\emptyset)\neq 0$.
Let ${\frak S}^l$ be the set of minimal length 
representatives of the cosets in $\widehat{\frak S}_l/{\frak S}_l.$ 
Then Lemma 9.3 implies that $m=yt$ with
$t\in{\frak S}_l$ and $y\in{\frak S}^\jb\cap{\frak S}^l$ such that
$S_\jb y\cap y S_l=\emptyset.$ Then,
$$\fb_{O_\nb}(\wedge\x_\emptyset)=
v^{l(\omega_\jb)}e_\jb\tilde T_y\tilde T_t\tilde T_\omega e^-=
v^{l(\omega_\jb)}e_\jb\tilde T_y\tilde T_\omega\tilde T_{\omega t\omega}e^-=
(-v)^{l(t)}\wedge\x_{(\jb)y\omega}\in\bigoplus_\lambda\SS\,|\lambda\rangle.$$
As a consequence if $l\leq n$ then $\Bb_l^+=\Bb_l$.

\vskip3mm

\noindent{\bf Conjecture.} {\it The bases $\Bb_l$ and $\Bb_l^+$ coincide
for all $l$.}\qed

\vskip5mm

\centerline{\bf 10. Proof of Theorem 6.3.}

\vskip3mm

\noindent{\bf 10.1.}
Let $\bigotimes^\infty$ be the free $\AA$-module linearly generated
by the semi-infinite monomials
$$\otimes\x_\ib=\x_{i_1}\otimes\x_{i_2}\otimes\x_{i_3}\otimes\cdots$$
where $\ib=(i_1,i_2,...)$ is a sequence of integers such that $i_k=1-k$
for $k>>1$. The affine Hecke algebra of type ${\frak{gl}}_\infty$
acts on $\bigotimes^\infty$ via formulas $(8.b)$ and $(8.c)$.
Set $\Omega^\infty=\sum_i\Im(T_i+1)\subset\bigotimes^\infty$.
As above $\wedge\x_\ib$ is the class of $\otimes\x_\ib$ in
$\bigotimes^\infty/\Omega^\infty$.
The formulas in Section 8 and [KMS, Lemma 2.2] imply that for all 
$\alpha\in\NN^{\ZZ/n\ZZ}$ we have
$$\forall\ib\quad\exists l\in\NN^\times\quad\text{such\ that}\quad
\fb_{\alpha}(\wedge\x_\ib)=
\fb_{\alpha}(\x_{i_1}\wedge\cdots\wedge\x_{i_l})\wedge\x_{i_{l+1}}
\wedge\x_{i_{l+2}}\wedge...\leqno{(a)}$$
Thus the action of $\Ub^-_n$ on $\bigwedge^l$ induces an action on $\fock$.

\vskip3mm

\noindent{\bf Lemma.} {\it The map 
$|\lambda\rangle\mapsto\wedge\x_\ib$, where $i_k=1+\lambda_k-k$,
gives an embedding of the representation of $\Ub^-_n$
on $\fock$ given in Section 6 into $\bigotimes^\infty/\Omega^\infty$.}

\vskip3mm

\noindent{\it Proof.} The proof goes by a direct computation.
First observe that $\wedge\x_\ib$ and $|\lambda\rangle$ have the same weight
for any $\lambda\in\Pi$ if $\ib$ is the sequence such that $i_k=1+\lambda_k-k$.
Fix $\lambda\in\Pi$ and $\alpha\in\NN^{\ZZ/n\ZZ}$. 
Fix $\ib$ as above. Formula $6.2.a$ gives
$$\fb_\alpha(\wedge\x_\ib)=
\fb_\alpha(|\lambda\rangle)=\sum_{a\,s.t.\,\ao=\alpha}
\gamma_a(\fb_\alpha)\prod_{j<i\atop\io=\jo}\kb_i^{a_j}
(|\lambda\rangle).$$ 
Moreover Remark 6.1 implies that $\gamma_a(\fb_\alpha)$ is $v^{h(a)}$ times the
product of the $\fb_i^{(a_i)}$'s ordered from $i=-\infty$ to $\infty$.
Using the formulas in Section 4 we first observe that $\fb_i^{(2)}$
acts by zero on the Fock space for any $i$.
The elements $|\lambda\rangle$ and $\otimes\x_\ib$ have the
same weight with respect to $\kb_i$. Thus we get
$$\fb_\alpha(\wedge\x_\ib)=\sum_\nb v^{e(\ib,\ib+\nb)}\wedge\x_{\ib+\nb},$$ 
where $\nb=(n_1,n_2,...)\in\{0,1\}^{\NN^\times}$ describes the set of all 
sequences such that $\alpha=\sum_{s\geq 1}n_s\epsilon_{\io_s}$ and
$$e(\ib,\ib+\nb)=\sum_{i_k>i_l}n_l\delta_{\io_l,\io_k}-
\sum_{i_k>1+i_l}n_l\delta_{\io_l+1,\io_k}-$$
$$-\sum_{i_k<i_l}n_ln_k\delta_{\io_l,\io_k}+
\sum_{1+i_k<i_l}n_ln_k\delta_{\io_l,\io_k+1}.$$
If $\wedge\x_{\ib+\nb}\neq 0$ then $e(\ib,\ib+\nb)=\sum_{i_k>i_l}n_l(1-n_k)
(\delta_{\io_l,\io_k}-\delta_{\io_l+1,\io_k})$.
On the other hand the formula in Section 8.3 gives
$$\fb_\alpha(\wedge\x_\ib)=\sum_\nb v^{c(\ib,\ib+\nb)}\wedge\x_{\ib+\nb},$$ 
where $\nb$ describes the same set and 
$$c(\ib,\ib+\nb)=-\sum_{1\leq k<l}n_l(1-n_k)
n(\epsilon_{\io_l},\epsilon_{\io_k}).$$ 
We are through (recall that $\ib$ is decreasing).
\qed

\vskip3mm

\noindent Theorem 6.3 follows from Proposition 9.2 and Lemma 10.1.

\vskip3mm

\noindent{\bf 10.2.} 
The involution $\psi$ on $\bigwedge^l$ induces the semilinear 
involution $\psi$ on $\fock$ such that,
$${\ts\forall\ib,\quad l\geq\sum_k(i_k-1+k)\Rightarrow\psi(\wedge\x_\ib)=
\psi(\x_{i_1}\wedge\cdots\wedge\x_{i_l})\wedge\x_{i_{l+1}}\wedge\x_{i_{l+2}}
\wedge ...}$$
Proposition 9.3 implies that $\psi$ coincides with the involution 
on $\fock$ used in [LT].
In [LT] Leclerc and Thibon have defined two bases
$\Bb^\pm=\{\bb^\pm_\lambda\,|\,\lambda\in\Pi\}$ 
in $\fock$ such that for all $\lambda$
$$\psi(\bb^\pm_\lambda)=\bb^\pm_\lambda,\qquad
\bb^-_\lambda-|\lambda\rangle\in\bigoplus_{\mu<\lambda}
v^{-1}\bar\SS |\mu\rangle
\quad\text{and}\quad
\bb^+_\lambda-|\lambda\rangle\in\bigoplus_{\mu<\lambda}
v\,\SS |\mu\rangle.$$

\vskip2mm

\noindent Thus, Conjecture 9.5 is equivalent to 

\vskip3mm

\noindent{\bf Conjecture.} {\it The bases $\Bb$ and $\Bb^+$ coincide.}
\qed

\vskip3mm

\noindent Set
$\bb_\lambda=\sum_\mu d_{\mu\lambda}|\mu\rangle$ and
$\bb^\pm_\lambda=\sum_\mu e^\pm_{\mu\lambda}|\mu\rangle$.
Conjecture 10.2 is precisely $d_{\lambda \mu}=e^+_{\lambda\mu}.$
\vskip5mm

\centerline{\bf 11. Proof of the Decomposition Conjecture.} 

\vskip3mm

\noindent Let $\varepsilon$ be a $n$-th root of unity.
The quantized Schur algebra $\Sb_{n,l}$ is the subalgebra
of $\widehat\Sb_{n,l}$ spanned by the elements $T_\mb$ with 
$\mb\in{\frak S}_\ib\setminus{\frak S}_l/{\frak S}_\jb$
(see Subsection 7.4). Fix $l\leq k$. Let consider the subalgebra
$$\Sb_l=\Sb_{k,l}\cap\bigoplus_{\lambda,\mu\in\Pi(l)}\widehat
\Hb_{\ib_\lambda\ib_\mu}$$
where $\ib_\lambda=(1-k)^{\lambda_1}(2-k)^{\lambda_2}...0^{\lambda_k}$
for any $\lambda=(\lambda_1\geq\lambda_2\geq\ldots)\in\Pi(l).$ 
We want to compute the decomposition matrices of the simple
$\Sb_{l}$-modules under the specialization $v\mapsto\varepsilon.$
The algebra $\Sb_{k,l}$ is Morita equivalent to $\Sb_l.$ 
For any $t\in\CC^\times$ let $\Sb_{l|t}$ be the specialization of 
$\Sb_{l}$ at $v=t$. 
The simple modules of $\Sb_{l|t}$ are parametrized by $\Pi(l).$
For any $k$ let $\Ub({\frak{gl}}_k)$ be the Lusztig integral form
of the quantized enveloping algebra of ${\frak{gl}}_k$ and let
$\Ub_\var({\frak{gl}}_k)$ be the specialization at $v=\var$. 
The set $\Pi_k$ is identified with the set of dominant weights of 
${\frak{gl}}_k$ with non-negative components.
If $\lambda\in\Pi_k$, let $V_\lambda$ and $W_\lambda$ be respectively
the simple and the Weyl $\Ub_\var({\frak{gl}}_k)$-module with 
highest weight $\lambda$. There exists a surjective map 
$\pi:\Ub_\var({\frak{gl}}_k)\to\Sb_{k,l|\var}$ (see [D2]).
If $\lambda\in\Pi(l)$ let $L_\lambda$, $M_\lambda$, be 
the simple and the Specht $\Sb_{k,l|\var}$-modules such that
$$\pi^*[L_\lambda]=[V_{\lambda'}]\quad\text{and}\quad 
\pi^*[M_\lambda]=[W_{\lambda'}]$$
in the Grothendieck ring.

\vskip3mm

\noindent{\bf Theorem.} {\it The specialization at $v=1$ of the matrix 
$(e^+_{\lambda\mu})_{\lambda\mu}$, $\lambda,\mu\in\Pi(l)$,
is the decomposition matrix of the Specht modules of $\Sb_l$.}

\vskip3mm

\noindent{\it Proof.} 
The Lusztig conjecture (proved by Kashiwara-Tanisaki and 
Kazhdan-Lusztig) gives the multiplicity of $W_\mu$
in $V_\lambda$. More precisely we have
$$[V_\lambda:W_\mu]=\sum_y(-1)^{l(yx)}P_{yx}(1),$$
where $x\in\widehat{\frak S}_l$ is minimal such that $\nu=\lambda\cdot x^{-1}$
satisfies
$$\nu_i<\nu_{i+1}\quad\forall i=1,2,...,k-1,\qquad\nu_1-\nu_k\geq 1-k-n,$$
and $\mu=\lambda\cdot x^{-1}y$.
According to Theorem $9.5.a$, the Lusztig Conjecture is equivalent to 
$$[L_{\lambda'}]=\sum_\mu e^-_{\lambda\mu}(1)\,[M_{\mu'}],\qquad
\forall\lambda\in\Pi_k.\leqno(a)$$ 
Recall that
$(e^+_{\lambda\mu})_{\lambda\mu}=(\bar e^-_{\lambda'\mu'})^{-1}_{\lambda\mu}$
(see [LT, Section 4]). Thus 
$$(a)\iff [M_{\lambda}]=\sum_\mu e^+_{\lambda\mu}(1)\,[L_\mu].$$\qed

\vskip5mm

\centerline{\bf 12. The Lusztig conjecture.}

\vskip3mm

\noindent{\bf 12.1.} 
Let $F$ be the variety of partial flags in $\CC^l$ of the type
$$\{0\}\subseteq F_1\subseteq F_2\subseteq\cdots\subseteq F_k=\CC^l.$$
The linear group $GL_l$ acts diagonaly on $F\times F$. Let 
$Z\subset T^*F\times T^*F$ be the Steinberg variety ($Z$ is a reducible
variety whose irreducible components are the closure of the conormal
bundles to the $GL_l$-orbits in $F\times F$).
The group $G=GL_l\times\CC^\times$ acts naturally on $Z$ :
the linear group acts diagonaly and $z\in\CC^\times$ acts by multiplication
by $z^{-2}$ along the fibers. The complexified Grothendieck group
of equivariant coherent sheaves on $Z$, denoted by $\Kb_{k,l}$, is endowed
with an associative convolution product (see [GV], [V2]) denoted by $\star$. 
For any $z\in\CC^\times$, a 
parametrization of the simple modules of the specialized algebra
$\Kb_{k,l|v=z}$ is given in [GV] (see in [V2] the remark after Theorem 4
for the case of roots of unity): 
the simple modules are labelled by orbits of
pairs $(s,x)\in\GL_l\times{\frak{gl}}_l$ 
where $s$ is semi-simple, $x^k=0$, and $sxs^{-1}=z^{-2}x$. 
As usual the $GL_l$-orbit of $x$ is labelled by the partition 
$\lambda\in\Pi(l)$ such that $\lambda_i$ is the length of the $i$-th Jordan
block of $x$. Then $\lambda'\in\Pi(l)\cap\Pi_k$.
The orbits of the pairs $(s,x)$ such that the spectrum of $s$ is in $z^{2\ZZ}$
are labelled by isomorphism class of nilpotent representations 
of $\Gamma_\infty$ if $z$ is generic and of $\Gamma_n$ if $z=\var$
(recall that $\var^2$ is a primitive $n$-th root of unity). 
Let $\Omega_{k,l}$ and $\Omega^\infty_{k,l}$ be the corresponding sets
of isomorphism classes of representations of $\Gamma_n$ and $\Gamma_\infty$.
If $O\in\Omega_{k,l}^\infty$ (resp. $O\in\Omega_{k,l}$) let
$L_O^\infty$ (resp. $L_O$) be the simple $\Kb_{k,l}$-module labelled by $O$. 
Similarly let $M_O^\infty$ and $M_O$ be the standard modules
labelled by $O$ (see [V2]).
Let $[M]$ be the class of the module $M$ in the complexified Grothendieck ring.
Let $\widehat\Rb_n$ and $\widehat\Rb_\infty$ be the linear span 
of the elements $[L_O]$ and $[L_O^\infty]$ where $O\in\Omega_{k,k}$
or $O\in\Omega^\infty_{k,k}$ and $k\geq 1$.
The restricted dual $\widehat\Rb_n^*$ (resp. $\widehat\Rb_\infty^*$)
is spanned by the linear forms $l_O$ (resp. $l^\infty_O$) such that
$$l_O([L_{O'}])=\delta_{O,O'}
\quad\text{and}\quad
l^\infty_O([L_{O'}^\infty])=\delta_{O,O'}.$$

\vskip3mm

\noindent{\bf 12.2.}
The quantized enveloping algebra of $\widehat{\frak{gl}}_k$ is generated
by elements $\eb_{i,s}$, $\fb_{i,s}$, $\hb_{j,t}$ and $\kb^{\pm 1}_j$
($0<i<k$, $0<j\leq k$, $s\in\ZZ$, $t\in\ZZ^\times$) which satisfy
the relations of the Drinfeld new presentation.
Let $\Ub(\widehat{\frak{gl}}_k)$ be the $\AA$-subalgebra generated
by the elements 
$\eb_{i,s}^{(m)}$, $\fb_{i,s}^{(m)}$, $[t]^{-1}\hb_{j,t}$ and $\kb^{\pm 1}_j$.
For any $z\in\CC^\times$ let $\Ub_z(\widehat{\frak{gl}}_k)$ be its 
specialization at $v=z$. 
In [GV], [V2], is defined a surjective algebra homomorphism
$\Psi_{k,l}\,:\,\Ub(\widehat{\frak{gl}}_k)\otimes_\AA\CC(v)\to 
\Kb_{k,l}\otimes_\AA\CC(v)$. It is proved in [S] that $\Psi_{k,l}$ 
restricts to a surjective homomorphism $\Ub(\widehat{\frak{gl}}_k)\to\Kb_{k,l}$.
Observe that the restriction of a simple $\Ub_z(\widehat{\frak{gl}}_k)$-module 
to $\Ub_z(\widehat{\frak{sl}}_k)$ is simple. Thus $\Psi^*_{k,l}L_O$ for
$O\in\Omega_{k,l}$ (resp. $\Psi^*_{k,l}L_O^\infty$ for $O\in
\Omega^\infty_{k,l}$), 
may be viewed as a simple $\Ub_z(\widehat{\frak{sl}}_k)$-module when $z=\var$ 
(resp. $z$ generic).
Recall that there is an algebra homomorphism
$ev\,:\,\Ub(\widehat{\frak{sl}}_k)\to\Ub({\frak{gl}}_k)$ such that
$$\matrix
ev(\eb_0)=v^{-1}\{\fb_{k-1},\{\fb_{k-2},...\{\fb_2,\fb_1\}...\}\}\kb_k\kb_{k-1} 
\hfill\cr\cr
ev(\fb_0)=(-1)^kv^{k-1}\{\eb_{k-1},\{\eb_{k-2},...\{\eb_2,\eb_1\}...\}\}
\kb_k^{-1}\kb_{k-1}^{-1} 
\hfill\cr\cr
ev(\fb_i)=\fb_i,\qquad ev(\eb_i)=\eb_i,\qquad i=1,2,...,k-1,\hfill
\endmatrix
$$
where $\{x,y\}=xy-v^{-1}yx$. If $\lambda\in\Pi_k$ let
$V_{\lambda}$ (resp. $V_{\lambda}^\infty$) be the simple
$\Ub_z({\frak{gl}}_k)$-modules with highest weight $\lambda$
where $z=\var$ (resp. $z$ generic).
The Drinfeld polynomials of $L_O$ and $L_O^\infty$ are computed in [V2].
If $\lambda\in\Pi(k),$ then
$\Psi^*_{k,k}L_{O_\lambda}$ and $\Psi^*_{k,k}L_{O_\lambda}^\infty$ are the 
pull-backs of the modules $V_{\lambda'}$ and $V_{\lambda'}^\infty$ 
by the evaluation map $ev$ (see [CP, Proposition 12.2.13]).
Let $\Rb_n$ and $\Rb_\infty$ be the linear span of the classes
$[V_\lambda]$ and $[V_\lambda^\infty]$ for all $\lambda$ and all $k$.
The restricted dual spaces $\Rb_n^*$ and $\Rb_\infty^*$ are spanned
by the linear forms $l_\lambda$ and $l^\infty_\lambda$ such that
$$l_\lambda([V_{\mu'}])=\delta_{\lambda\mu}
\quad\text{and}\quad
l^\infty_\lambda([V_{\mu'}^\infty])=\delta_{\lambda\mu}.$$
The element $[V_\lambda^\infty]$ may be viewed as the class in $\Rb_n$ of
the Weyl module $W_\lambda$ with highest weight $\lambda$. 
Let $s_\lambda\in\Rb_n^*$ be such that
$$s_\lambda([W_{\mu'}])=\delta_{\lambda\mu}.$$

\vskip3mm

\noindent{\bf 12.3.}
In this subsection $\Ub^-_n$, $\Ub^-_\infty$ and $\fock$ 
stand for their specializations at $v=1$.

\vskip3mm

\noindent{\bf Theorem.} {\it The linear isomorphism 
$\Rb_n^*{\buildrel\sim\over\to}\fock$
such that $s_\lambda\mapsto|\lambda\rangle$ maps
$l_{\lambda}$ to $\bb_\lambda$.}

\vskip3mm

\noindent{\it Proof.}
First observe that the classes of the standard modules $[M_O]$ and 
$[M_O^\infty]$ form a basis of the spaces $\widehat\Rb_n$ and 
$\widehat\Rb_\infty$. Let $m_O$ and $m_O^\infty$ be the elements
of the dual basis.
To avoid confusions let $\fb_O^\infty$, $\bb_O^\infty$, denote the 
generatrors of $\Ub_\infty^-$.
The multiplicity formula [GV, Theorem 6.6] implies that there are
two linear isomorphisms 
$$\iota_n\,:\,\Ub_{n}^-\to\widehat\Rb_n^*
\quad\text{and}\quad
\iota_\infty\,:\,\Ub_{\infty}^-\to\widehat\Rb_\infty^*$$
such that
$$\iota_n(\fb_O)=m_O,\quad\iota_n(\bb_O)=l_O,\quad
\iota_\infty(\fb^\infty_O)=m^\infty_O,\quad
\iota_\infty(\fb^\infty_O)=l^\infty_O.\leqno(a)$$
The spaces $\Rb_n^*$ and $\Rb_\infty^*$ are identified with
$\bigwedge^\infty$ via the maps  
$$s_{\lambda}\mapsto|\lambda\rangle
\quad\text{and}\quad
l^\infty_{\lambda}\mapsto|\lambda\rangle.$$
We obtain the following commutative square
$$\matrix
\bigwedge^\infty&=&\Rb_n^*&{\buildrel\sim\over\to}&\Rb_\infty^*&
=&\bigwedge^\infty\cr
&&\uparrow&&\uparrow&&\cr
\Ub^-_{n}&=&\widehat\Rb_n^*&\to&\widehat\Rb_\infty^*&=&\Ub^-_{\infty},
\endmatrix$$
where the horizontal arrows are the dual of the specialization maps  
and the vertical arrows are the dual of the evaluation maps.
By definition the upper arrow maps $s_\lambda$ to $l^\infty_\lambda$ and 
both elements are identified with the vacuum vector $|\lambda\rangle$. 
The right vertical arrow is such that 
$$\bb^\infty_O=l^\infty_O\mapsto l^\infty_{\lambda}=|\lambda\rangle
\quad\text{if}\quad O=O_\lambda,\qquad
\bb_O^\infty\mapsto 0\quad\text{else}.$$
By Proposition 5 it is the quotient map
$${\ts\Ub^-_{\infty}\to\bigwedge^\infty,
\qquad u\mapsto u(|\emptyset\rangle).}$$
Suppose first that the lower horizontal arrow is the map
$\gamma$ introduced in Section 6. Then the left 
vertical arrow is the quotient map
$${\ts\Ub^-_{n}\to\bigwedge^\infty,\qquad u\mapsto u(|\emptyset\rangle)}.$$
Hence $(a)$ implies that the left vertical arrow maps $l_{O_{\lambda}}$ to 
$\bb_\lambda$. Since this arrow is the transpose of the evaluation map we get
$l_{\lambda}=\bb_\lambda$ and we are through. 
By Subsection 6.4, to prove that the map 
$\widehat\Rb^*_n\to\widehat\Rb_\infty^*$ is $\gamma$ 
we are reduced to prove that if $r(O')\subseteq O$ then
$[M^\infty_{O'}]$ specializes to $[M_O]$. This is obvious by
the localization theorem in equivariant $K$-theory.
\qed

\vskip 3mm

\noindent{\bf 12.4.}
Theorem 12.3 implies that 
$[V_{\lambda}^\infty]=\sum_\mu d_{\lambda'\mu'}(1)\,[V_\mu].$ 
According to Section 11 the Lusztig Conjecture can be written as 
$$[W_{\lambda}]=\sum_\mu e^+_{\lambda'\mu'}(1)\,[V_\mu]\quad
\forall\lambda,$$ 
which is precisely  Conjecture 10.2.

\vskip5mm

\centerline{\bf 13. Proof of Proposition 6.1.}

\vskip3mm

\noindent{\bf 13.1.}
Fix $\Gamma=\Gamma_n$ or $\Gamma_\infty$.
Let $\Sc_d$ be the set of finite sequences 
$\db=(d^1,d^2,...,d^l)$ of elements in $\NN^{(I)}$ 
such that $\sum_kd^k=d$.
Fix a $I$-graded vector space $V$ of dimension $d$. 
For each $\db\in\Sc_d$ let $F_\db$ be the set of flags of 
$V$ of type $\db$, i.e. $F_\db$ is the set of filtrations 
$F=(\{0\}=F^0\subseteq F^1\subseteq\cdots\subseteq F^l=V)$
such that $F^k$ is $I$-graded and has dimension $d^1+d^2+\cdots+d^k$.
Given $x\in E_V$ we say that a flag $F\in F_\db$ is $x$-stable 
if $x(F^k)\subseteq F^{k-1}$ for all $k$. 
Let $\tilde F_\db$ be the variety of all pairs $(x,F)$ 
such that $x\in E_V$ and $F\in F_\db$ is $x$-stable.
The group $G_V$ acts on $\tilde F_\db$ in the obvious way. 
Let $\pi_\db\,:\,\tilde F_\db\to E_V$ be the first projection. 
The map $\pi_\db$ commutes to $G_V$.
Thus the function $f_\db=\pi_{\db\,!}(1)$ belongs to $\CC_{G_V}(E_V)$.

\vskip3mm

\noindent{\bf Lemma.} {\it 
$(a)$ The space $\CC_{G_V}(E_V)$ is linearly spanned
by the elements $f_\db$ with $\db\in\Sc_d$.\hfill\break
$(b)$ For any $a,b\in\NN^{(I)}$ and any $\ab\in\Sc_a$,
$\bb\in\Sc_b$, we have $f_\ab\circ f_\bb=q^{-m(b,a)}f_{\ab\bb}$ where
$\ab\bb\in\Sc_{a+b}$ is the sequence $\ab$ followed by the sequence $\bb$. 
}

\vskip3mm

\noindent{\it Proof.} Claim $(b)$ is proved as in [L2, Lemma $3.2.b$].
Let us prove claim $(a)$. 
If a flag $F$ is $x$-stable then $F^k\subseteq\Ker(x^k)$. Thus if
$\db\in\Sc_d$ is such that 
$$d^1+d^2+\cdots+d^k=\dim\Ker(x^k)\quad\forall k=1,2,3,...,$$
then $\pi_\db^{-1}(x)$ is reduced to the single flag
$$\{0\}\subseteq\Ker(x)\subseteq\Ker(x^2)\subseteq\cdots\subseteq V.$$
In particular $f_\db(x)=1$. 
Moreover, in this case $f_\db$ is supported on the
$G_V$-orbits of the $y$'s such that
$$\dim\Ker(x^k)\leq\dim\Ker(y^k)\quad\forall k=1,2,3,...,$$
i.e. $y\in\overline{G_V\cdot x}$. We are through.
\qed

\vskip3mm

\noindent{\bf Remark.} It is easy to see that for any $d\in\NN^{\ZZ/n\ZZ}$
and $\db=(d)$ we have $f_\db$=$\fb_d$. Thus Proposition 3.5 is a consequence of
$(a)$ and $(b)$.

\vskip3mm

\noindent{\bf 13.2.} 
We fix a $\ZZ$-graded vector space $V$ of dimension $d$.
Let $\Vo$ be the associated $\ZZ/n\ZZ$-graded vector space, of dimension $\do$.
The space $\Vo$ is endowed with the $\ZZ$-filtration whose associated graded
is identified with $V$. Fix $\dbo\in\Sc_\do$. 
We have the following commutative diagram 
$$\matrix
\tilde F_\dbo&{\buildrel\pi_\dbo\over\lra}&E_\Vo&&\cr
\uparrow&&\uparrow{\scriptstyle j}&&\cr
\tilde F_{\dbo,d}&{\buildrel\pi_{\dbo,d}\over\lra}&E_{\Vo,V}&
{\buildrel p\over\lra}&E_V,
\endmatrix$$
where $\tilde F_{\dbo,d}=\pi_\dbo^{-1}(E_{\Vo,V})$ 
and the vertical arrows are the embeddings. We have clearly
$$p_!j^*(f_\dbo)=(p\pi_{\dbo,d})_!(1).\leqno(c)$$
Let $\Sc_{\dbo,d}\subset\Sc_d$ be the set of sequences $\db$ such that
$\sum_{i\in\io}d^k_i=\do^k_\io$ for any $k$ and $\io$.
If $\db\in\Sc_{\dbo,d}$ let $\tilde F_{\dbo,\db}\subset\tilde F_{\dbo,d}$ 
be the set of pairs $(x,F)$ such that the associated graded 
of $F^k$ with respect to the filtration induced by the 
$\ZZ$-filtration on $\Vo$ has dimension $d^1+d^2+\cdots+d^k$.
The sets $\tilde F_{\dbo,\db}$ form a partition of $\tilde F_{\dbo,d}$. 
We have a commutative square
$$\matrix
\tilde F_{\dbo,d}&{\buildrel p\pi_{\dbo,d}\over\lra}&E_V\cr
\uparrow&&\uparrow{\scriptstyle\pi_\db}\cr
\tilde F_{\dbo,\db}&{\buildrel\tau\over\lra}&\tilde F_\db,
\endmatrix$$
where the left vertical arrow is the inclusion and $\tau$ maps the pair 
$(x,F)$ to the associated graded. Thus
$$(p\pi_{\dbo,d})_!(1)=\sum_{\db\in\Sc_{\dbo,d}}(\pi_\db\tau)_!(1).\leqno(d)$$

\noindent{\bf Lemma.} {\it The map $\tau$ is a vector bundle of rank
$$r(\db)=\sum_{k\geq l}\sum_{i>j\atop\io=\jo}d_j^kd_{i+1}^l+
\sum_{k<l}\sum_{i>j\atop\io=\jo}d_j^kd_i^l.$$}

\noindent{\it Proof.} The proof goes as the proof of [L2, Lemma 4.4].
More precisely fix $(x,F)\in\tilde F_\db$ and 
compute the fiber $\tau^{-1}(x,F)$.
Giving a $\ZZ/n\ZZ$-graded subspace 
$\bar F^k\in\bar V$ of dimension $\do^1+\do^2+\cdots+\do^k$ whose 
associated $\ZZ$-graded is $F^k$ is the same as giving a map 
$$z^k=\oplus z^k_{ij}\in\bigoplus_{i>j\atop\io=\jo}\Hom(F^k_j,V_i/F^k_i).$$
Then $\bar F^k\subset\bar F^{k+1}$ if and only if 
$z^{k+1}=z^k\,:\,F^k\to V/F^{k+1}.$
On the other hand giving $\bar x\in E_{\bar V,V}$ such that 
$p(\bar x)=x$ is the same as giving a map
$$y=\oplus y_{i+1,j}\in\bigoplus_{i>j\atop\io=\jo}\Hom(V_j,V_{i+1}).$$
Then $\bar F$ is $\bar x$-stable if and only if 
$$z^k_{i+1,j+1}\circ x_j-x_i\circ z^k_{ij}-y_{i+1,j}=0\,:\,
F^k_j\to V_{i+1}/F^k_{i+1}.$$
The Lemma results from a direct computation.
\qed

\vskip3mm

\noindent The Lemma and $(c)$, $(d)$, give
$${\ts\gamma_d(f_\dbo)=\sum_{\db\in\Sc_{\dbo,d}}q^{2r(\db)-h(d)}f_\db.}$$
Fix $\alpha,\beta\in\NN^{\ZZ/n\ZZ}$, $\abo\in\Sc_\alpha$, and
$\bbo\in\Sc_\beta$. Using Lemma $13.1.b$ we get
$${\ts\gamma_d(f_\abo\circ f_\bbo)=
\sum_{\ab,\bb}q^{m(b,a)-m(\beta,\alpha)+2r(\ab\bb)-h(d)}f_\ab\circ f_\bb,}$$
where the sum is over all $(\ab,\bb)\in\Sc_{\abo,a}\times\Sc_{\bbo,b}$
and all $(a,b)$ such that $\ao=\alpha,$ $\bo=\beta$, and $d=a+b$.
We are thus reduced to prove the following identity
$$m(b,a)-m(\beta,\alpha)+2r(\ab\bb)-2r(\ab)-2r(\bb)+h(a)+h(b)-h(d)
=k(b,a).\leqno(e)$$
Set
$$l_+(b,a)=\sum_{i>j\atop\io=\jo}(b_ia_j+b_ja_{i+1})
\quad\text{and}\quad
l_-(b,a)=\sum_{i<j\atop\io=\jo}(b_ia_j+b_ja_{i+1}).$$
Then $(e)$ follows from the following equalities which are easy to
prove :
$$\matrix
m(b,a)-m(\beta,\alpha)=-l_+(b,a)-l_-(b,a),\hfill\cr\cr
h(a)+h(b)-h(d)=k(b,a)-l_+(b,a)+l_-(b,a),\hfill\cr\cr
r(\ab\bb)-r(\ab)-r(\bb)=l_+(b,a).\hfill
\endmatrix$$

\vskip1cm

\noindent{\it Acknowledgements.}
{\eightpoint{This work was partially done while the second author
was visiting the University of Chicago. The second author is grateful 
to the departement of Mathematics of the U.C., in particular to V. Ginzburg, 
for his kind invitation.}}

\vskip1cm

\centerline{\bf References}

\vskip3mm

\hangindent3cm[A]\quad\qquad\  Ariki, S.: On the decomposition numbers of the
Hecke algebra of $G(m,1,n)$. {\sl J. Math. Kyoto Univ}, {\bf 36} (1996), 789-808.

\vskip1mm

\hangindent3cm[CG]\quad\qquad Chriss, N., Ginzburg, V.: 
Representation theory and complex geometry.
{\sl Birkhauser} (1997).

\vskip1mm

\hangindent3cm[CP]\quad\qquad Chari, V., Pressley, A.: 
A guide to quantum groups.
{\sl Cambridge University Press} (1994).

\vskip1mm

\hangindent3cm[D1]\quad\qquad\  Deodhar, V.V.: On some geometric aspects of the 
Bruhat ordering II. The parabolic analogue of Kazhdan-Lusztig polynomials. 
{\sl J. Algebra}, {\bf 111} (1987), 483-506.

\vskip1mm

\hangindent3cm[D2]\quad\qquad\  Du, J.: A note on quantized Weyl
reciprocity at roots of unity.
{\sl Algebra Colloq.}, {\bf 2} (1995), 363-372.

\vskip1mm

\hangindent3cm[G1]\quad\qquad Green, J.A.: 
Hall algebras, hereditary algebras and quantum groups.
{\sl Invent. Math.}, {\bf 120} (1995), 361-377.

\vskip1mm

\hangindent3cm[G2]\quad\qquad Green, R.M.: 
The affine $q$-Schur algebra.
{\sl qalg-preprint}, {\bf 9705015}.

\vskip2mm

\hangindent3cm[GV]\quad\qquad Ginzburg, V., Vasserot, E.: 
Langlands reciprocity for affine quantum groups of type $A_n$.
{\sl Internat. Math. Res. Notices}, {\bf 3} (1993), 67-85.

\vskip1mm

\hangindent3cm[H]\quad\qquad\  Hayashi, T.: $Q$-analogues of Clifford and 
Weyl algebras - spinor and oscillator representations of quantum
enveloping algebras. {\sl Comm. Math. Phys.}, {\bf 127} (1990), 129-144.

\vskip1mm

\hangindent3cm[IM]\quad\qquad\ Iwahori, N., Matsumoto, H.: 
On some Bruhat decompositions and the structure of Hecke rings of $p$-adic 
Chevalley groups. {\sl Pub. I.H.E.S.}, {\bf 25} (1965), 5-48.

\vskip1mm

\hangindent3cm[KMS]\quad\quad Kashiwara, M., Miwa, T., Stern, E.: Decomposition
of $q$-deformed Fock space. {\sl Selecta Mathematica, New Series}, {\bf 1}
(1995), 787-805.

\vskip1mm

\hangindent3cm[L1]\quad\qquad\ Lusztig, G.: Canonical bases arising from
quantized enveloping algebras. {\sl J. Amer. Math. Soc.}, {\bf 3} (1990), 
447-498. 

\vskip1mm

\hangindent3cm[L2]\quad\qquad\ Lusztig, G.: Quivers, perverse sheaves and 
enveloping algebras. {\sl J. Amer. Math. Soc.}, {\bf 4} (1991), 365-421. 

\vskip1mm

\hangindent3cm[L3]\quad\qquad\ Lusztig, G.: Introduction to quantum groups.
{\sl Birkh\"auser, Progr. in Math.}, {\bf 110} (1993).

\vskip1mm

\hangindent3cm[L4]\quad\qquad\ Lusztig, G.: Canonical bases and Hall algebras.
{\sl Preprint}, (1997).

\vskip2mm

\hangindent3cm[LT]\quad\qquad\ Leclerc, B., Thibon, J.-Y.: Canonical bases of
$q$-deformed Fock spaces. {\sl Internat. Math. Res. Notices}, {\bf 9} (1996),
447-456.

\vskip1mm

\hangindent3cm[MM]\quad\qquad Misra, K.C., Miwa, T.: Crystal base for the basic
representation of $U_q(\widehat{\frak{sl}}_n)$.
{\sl Comm. Math. Phys.}, {\bf 134} (1990), 79-88.

\vskip1mm

\hangindent3cm[N]\quad\qquad\ Nakajima, I.: Instantons on ALE spaces, quiver
varieties, and Kac-Moody algebras. {\sl Duke Math. J.}, {\bf 76} (1994),
365-416.

\vskip1mm

\hangindent3cm[S]\quad\qquad\  Schiffmann, O.: 
Alg\`ebres affines quantiques aux racines de
l'unit\'e et $K$-th\'eorie \'equivariante. 
{\sl C. R. Acad. Sci. Paris}, to appear.

\vskip1mm

\hangindent3cm[V1]\quad\qquad\ Vasserot, E.: Repr\'esentations de groupes
quantiques et permutations. {\sl Ann. Scient. \'Ec. Norm. Sup.}, 
{\bf 4\`eme s\'erie, 26} (1993), 747-773.

\vskip1mm

\hangindent3cm[V2]\quad\qquad\ Vasserot, E.: Affine quantum groups and 
equivariant $K$-theory. {\sl Transformation Groups}, to appear.

\vskip1mm

\hangindent3cm[VV]\quad\qquad Varagnolo, M., Vasserot, E.:
Double-loop algebras and the Fock space.	
{\sl Invent. Math.}, {\bf 133} (1998), 133-159.

\vskip1mm

\vskip3cm
{\eightpoint{
$$\matrix\format\l&\l&\l&\l\\
\phantom{.} & {\text{Michela Varagnolo}}\phantom{xxxxxxxxxxxxx} &
{\text{Eric Vasserot}}\\
\phantom{.}&{\text{D\'epartement de Math\'ematiques}}\phantom{xxxxxxxxxxxxx} &
{\text{D\'epartement de Math\'ematiques}}\\
\phantom{.}&{\text{Universit\'e de Cergy-Pontoise}}\phantom{xxxxxxxxxxxxx} &
{\text{Universit\'e de Cergy-Pontoise}}\\
\phantom{.}&{\text{2 Av. A. Chauvin}}\phantom{xxxxxxxxxxxxx} & 
{\text{2 Av. A. Chauvin}}\\
\phantom{.}&{\text{95302 Cergy-Pontoise Cedex}}\phantom{xxxxxxxxxxxxx} & 
{\text{95302 Cergy-Pontoise Cedex}}\\
\phantom{.}&{\text{France}}\phantom{xxxxxxxxxxxxx} & 
{\roman{France}}\\
&{\text{email: varagnol\@math.pst.u-cergy.fr}}\phantom{xxxxxxxxxxxxx} &
{\text{email: vasserot\@math.pst.u-cergy.fr}}
\endmatrix$$
}}
\enddocument